\documentclass[12pt]{amsart}
\usepackage[utf8]{inputenc}
\usepackage[T1]{fontenc}
\usepackage{amsfonts}
\usepackage[leqno]{amsmath}
\usepackage{amssymb, comment, float}
\usepackage[dvipsnames]{xcolor}
\usepackage[foot]{amsaddr}
\usepackage[shortlabels]{enumitem}
\usepackage{mathdots}
\usepackage{thm-restate}
\usepackage[footskip=0.5in, headheight = 0.5in, top=1.25in, bottom=1.25in,  right=1in,  left=1in]{geometry}
\usepackage{hyperref}
\hypersetup{
colorlinks,
linkcolor=black,
urlcolor=Blue,
citecolor=black,}
\usepackage{cleveref}
\usepackage[center]{caption}
\usepackage{eufrak}
\usepackage{marvosym}     
\usepackage{latexsym}
\usepackage[nodayofweek]{datetime}
\usepackage{tikz}
\usepackage{subfig}
\usepackage{thm-restate}
\usepackage{algorithmicx,algpseudocode}
\usepackage{mathtools}

\newtheorem{theorem}{Theorem}[section]
\newtheorem{lemma}[theorem]{Lemma}
\newtheorem{claim}[theorem]{Claim}
\newtheorem{problem}[theorem]{Problem}
\newtheorem{remark}[theorem]{Remark}

\newtheorem{corollary}[theorem]{Corollary}

\newtheorem{definition}[theorem]{Definition}
\newtheorem*{mainthmprime}{Theorem~\ref{thm:main} (primed notation)}

\newcommand{\poi}{\mathbb{N}}

\newtheorem{case}{Case}

\newcommand{\subcase}[1]{\vspace{4pt}{\bf Subcase {#1}}}

\newcounter{tbox}

\newcommand{\mylongtitle}[1]{%
  \ifodd\value{page}%
    \protect\parbox{0.97\linewidth}{#1}\hfill%
  \else%
    \hfill\protect\parbox{0.97\linewidth}{#1}%
  \fi%
}

\usepackage{graphicx}
\usepackage{eso-pic}

\renewcommand{\qedsymbol}{$\blacksquare$}

\title[Efficient Reconfiguration of Independent Transversals]{Efficient Reconfiguration of Independent Transversals}

\author{Sahab Hajebi}
\address{Department of Combinatorics and Optimization, University of Waterloo, Waterloo ON Canada N2L 3G1}
\email{s2hajebi@uwaterloo.ca}

\author{Penny Haxell}
\address{}
\email{pehaxell@uwaterloo.ca}
\thanks{Research of P. Haxell partially supported by NSERC}

\date{\today}

\begin{document}

\maketitle
\thispagestyle{plain}
\pagestyle{plain}



\begin{abstract}
We study the reconfiguration of independent transversals in graphs with a prescribed vertex partition. Given a graph $G$ with a vertex partition $\mathcal{V}=(V_1,\ldots,V_m)$, an \textit{independent transversal} is an independent set containing exactly one vertex from each part $V_i$ of $\mathcal{V}$. The reconfiguration problem concerns transforming one independent transversal into another through a sequence of local changes, while preserving the independent transversal property at every step. For any fixed $r\in\mathbb{N}$ and $\epsilon>0$, we present a polynomial-time algorithm that, given an $r$-claw-free graph $G$ with a vertex partition $\mathcal{V}=(V_1,\ldots,V_m)$ and two independent transversals $S$ and $R$, either finds a polynomial-length reconfiguration path between $S$ and $R$, or identifies a nonempty set of partition classes $\mathcal{B}\subseteq\mathcal{V}$ such that the subgraph induced by their union admits a strong dominating set of size at most $(2+\epsilon)|\mathcal{B}|$. Our result provides an algorithmic counterpart to recent connectivity results for the reconfiguration graph of all independent transversals, and establishes polynomial bounds on its diameter under suitable conditions. As a specific application, we obtain a polynomial-time algorithm for reconfiguring perfect matchings in $r$-uniform bipartite hypergraphs $H=(A,B,E)$, that either produces a polynomial-length reconfiguration path between two given perfect matchings or identifies a subset $S\subseteq A$ whose link hypergraph $E_S$ in $B$ has a small cover.
\end{abstract}

\section{Introduction}


Let $G$ be a graph with a partition $\mathcal{V}=(V_1, V_2,\dots,V_m)$ of $V(G)$. An \textit{independent transversal} (IT) $S$ for the pair $(G,\mathcal{V})$ is an independent set of $G$ that contains exactly one vertex from each part of $\mathcal{V}$. The notion of an IT provides a general framework that captures a wide range of combinatorial problems, since many such problems can be formulated by asking whether a given $(G,\mathcal{V})$ admits an independent transversal. We denote by $\mathcal{I}(G)$ the set of all independent sets of $G$, and by $\mathcal{I}(G,\mathcal{V})$ the set of all independent transversals for $(G,\mathcal{V)}$, respectively. We define the $\textit{reconfiguration graph}$ $\textbf{RG}(G, \mathcal{V})$ whose vertex set is $\mathcal{I}(G,\mathcal{V})$, in which two vertices $S,R$ are adjacent if and only if $|S\Delta R|=2$ and $R\cup S \in \mathcal{I}(G)$; that is, $S$ and $R$ differ in exactly one part $V_i\in \mathcal{V}$, and the corresponding vertices of $S$ and $R$ in $V_i$ are not adjacent in $G$. The general Reconfiguration problem asks whether one feasible solution can be transformed into another through a sequence of local changes while preserving feasibility at every step. In our setting, the question is whether every pair $(S,R)$ of independent transversals in $(G, \mathcal{V})$ is joined by a path in $\textbf{RG}(G, \mathcal{V})$. For two vertices $u,v\in V(G)$, a \emph{$u$--$v$ path} (\emph{walk}) is a path (walk) in $G$ whose endpoints are $u$ and $v$.

\begin{remark} Throughout this paper, we often use the term ``path'' for a walk between two vertices $u$ and $v$ without verifying that the walk has no repeated vertices. This suffices for our purposes, since we only seek an upper bound on the length of a shortest $u$--$v$ path, which we obtain by bounding the length of a $u$--$v$ walk. \end{remark}

Given a set of partition classes $\mathcal{B}\subseteq \mathcal{V}$, we let $ V_{\mathcal{B}}=\bigcup_{V_i\in \mathcal{B}}V_i. $  For vertex subsets $D,S\subseteq V(G)$, we say that $D$ \emph{strongly dominates} $S$ if every vertex of $S$ has a neighbor in $D$. Throughout this paper, the term \emph{domination} always refers to strong domination unless stated otherwise.  We denote by $\widetilde{\gamma}(S)$ the minimum size of a vertex set that strongly dominates $S$. 

A classical result of Haxell~\cite{Haxell_Match_HyperG, Haxell_Classic} guarantees the existence of an IT when the vertex classes are sufficiently large: if $G$ has maximum degree $\Delta$ and $|V_i|\geq 2\Delta$ for every $i$, then $(G,\mathcal{V})$ admits an IT. The algorithmic problem of finding independent transversals has also been extensively studied (see e.g.~\cite{Harris, HarrisSrinivasan}, and~\cite{GrafHaxell} for background on this topic). 
In particular, building on~\cite{Haxell_Match_HyperG, Haxell_Classic} and extending techniques introduced by Annamalai~\cite{Annamalai} for bipartite hypergraph matchings, Graf and Haxell~\cite{GrafHaxell} developed a polynomial-time algorithm for fixed $r$ and $\epsilon>0$ that, in $r$-claw-free graphs $(G,\mathcal{V})$, 
finds either an independent transversal, or a set of partition classes $\mathcal{B}\subseteq\mathcal{V}$ such that the subgraph $G[V_{\mathcal {B}}]$ of $G$ induced by  $V_{\mathcal{B}}$ satisfies $\tilde{\gamma} (G[V_{\mathcal {B}}])\leq(2+\epsilon)(|\mathcal{B}|-1)$. To say $(G,\mathcal{V})$ is $r$-\textit{claw-free} means that no vertex of $G$ has $r$ independent neighbors in distinct classes of the partition $\mathcal{V}$. Also, we denote $G[V_{\mathcal{B}}]$ by $G_{\mathcal{B}}$ for brevity. 

Buys, Kang, and Ozeki~\cite{BuyKangOzeki} established the following reconfiguration analogue of the consequence of~\cite{Haxell_Match_HyperG, Haxell_Classic} mentioned above.

\begin{theorem}\label{thmBKO}
Let $G$ be a graph of maximum degree $\Delta$, and let $\mathcal{V}=(V_1,\ldots,V_m)$ be a partition of $V(G)$ such that $|V_i|\geq 2\Delta$ for every $i\in[m]$. Suppose that, for every nonempty $I\subseteq[m]$, the subgraph $G[V_I]$ is not the disjoint union of $|I|$ copies of $K_{\Delta,\Delta}$. Then $\textbf{RG}(G, \mathcal{V})$ is connected.
\end{theorem}
In particular, their result implies that $\textbf{RG}(G, \mathcal{V})$ is connected whenever $|V_i|\geq 2\Delta+1$ for every $i\in[m]$. Their proof is combinatorial and does not provide an efficient algorithm for finding a reconfiguration path. They raise several natural quantitative questions about the reconfiguration graph, in particular, they ask for best possible upper bounds on the diameter of $\textbf{RG}(G, \mathcal{V})$ in terms of $|\mathcal V|=m$ and $\Delta$. This question is also closely related to one of the original motivations for studying reconfiguration. Buys, Kang, and Ozeki describe a natural Markov chain on independent transversals, in which one block is updated at a time. Their connectivity result guarantees the ergodicity of this chain, and they ask under what conditions it mixes in polynomial time.

The computational complexity of reconfiguration problems has been widely studied; see, for example, the survey of Nishimura~\cite{Nishimura}. In general, connectivity of a reconfiguration graph and a polynomial bound on its diameter are important initial steps, but they do not guarantee an efficient algorithm for constructing reconfiguration paths. 
Thus, obtaining efficient algorithms for constructing short reconfiguration paths is a natural strengthening of connectivity and diameter results.

The reconfiguration problem for independent transversals has been studied in various special cases, see e.g.~\cite{bonamy, cambie2024, cambie2025, cambie2026, wdowinski2026}. Wdowinski~\cite{wdowinski} recently established more general sufficient conditions, in the spirit of Theorem~\ref{thmBKO}, that guarantee connectivity of the corresponding reconfiguration graph. In particular, he proved the following theorem, where the parameter $i\gamma(S)$ is defined as the smallest integer $\ell$ such that every independent set of $G[S]$ is dominated in $G$ by a vertex set of size at most $\ell$.

\begin{theorem}\label{wdowinski}
   For a graph $G$ with a vertex partition $\mathcal{V}=(V_1,\ldots,V_m)$, the reconfiguration graph $\textbf{RG}(G,\mathcal{V})$ is connected if either of the following conditions holds: 
   \begin{itemize}
       \item [(a)] $\widetilde{\gamma}(V_I)\geq 2|I|+1,$ for every nonempty $I\subseteq[m]$;
       \item [(b)] $i\gamma(V_I)\geq |I|+1$ for every nonempty $I\subseteq[m]$.
   \end{itemize} 
\end{theorem}
The two conditions above arise naturally from the corresponding existence results for independent transversals: Haxell's theorem~\cite{Haxell_Match_HyperG, Haxell_Classic} in the first case, and the Aharoni--Haxell theorem~\cite{Haxell_HallThm} in the second. Although these conditions guarantee connectivity, again the proof is non-algorithmic and therefore does not directly provide an efficient procedure for constructing a reconfiguration path between two given independent transversals. Wdowinski also asked about the diameter of $\textbf{RG}(G,\mathcal{V})$, and in addition posed the problem of whether reconfiguration paths in it can be found efficiently. In this paper, we develop an algorithmic counterpart of Theorem~\ref{wdowinski}(a) for $r$-claw-free graphs, where $r$ is a fixed constant. Given two independent transversals, our algorithm runs in polynomial time and either constructs a polynomial-length path between them in the reconfiguration graph or produces a small dominating set for a subgraph induced by a set of partition classes. Some of the main ideas behind our algorithm are inspired by the framework of Graf and Haxell~\cite{GrafHaxell}, which we build upon and adapt to the reconfiguration setting.

Our main theorem can be stated as follows.

\begin{theorem}\label{thm:main}
    Let $r\in \poi$ and $\epsilon >0$ be given. Then there is an algorithm $A_{r,\epsilon}$ that takes an $r$-claw-free graph $G$ with respect to a vertex partition $\mathcal{V}=(V_1, V_2,\dots,V_m)$ and two independent transversals $S,R\in \mathcal{I}(G,\mathcal{V})$ as input, and finds in polynomial time either;
    \begin{itemize}
        \item [-] a polynomial-length $S$--$R$ path $P$ in $\textbf{RG}(G, \mathcal{V})$, or
            \item[-] a set of partition classes $\mathcal{B}\subseteq \mathcal{V}$ and a vertex set $D\subseteq V(G)$ such that $|D|\leq (2+\epsilon)|\mathcal{B}|$ and $D$ dominates $G_{\mathcal{B}}$. Furthermore, there exists a biconstellation $C$ for some $\mathcal{B}_0\supseteq \mathcal{B}$ such that $V(C)\subseteq D$ and $|D\setminus V(C)|< \epsilon(|\mathcal{B}|-1)$.
    \end{itemize}
\end{theorem}
Here, a \emph{biconstellation} is a technical structural notion, whose precise definition is given in Definition~\ref{biconstellation}. This additional structural information is useful in applications of our main theorem; in particular, we use it in Section~\ref{sec:Application} to obtain a polynomial-length reconfiguration path between two hypergraph matchings. 

At a high level, the algorithm $ A_{r,\epsilon}$ repeatedly tries to bring the two given independent transversals $S$ and $R$ closer together by increasing their common part. For the current pair of transversals, it builds an auxiliary structure called an \textit{alternating tree} to search for a sequence of valid reconfiguration steps that brings them closer together. If such a sequence is found, the algorithm updates the two transversals accordingly and repeats the process. Otherwise, the failure of this search implies that the current alternating tree, together with a small number of extra vertices (essentially forming the error term $\epsilon|\mathcal{B}|$),  provides a dominating set $D$ for the set of partition classes $\mathcal{B}$ that it intersects. These are the sets appearing in the second outcome of Theorem~\ref{thm:main}. Thus, the algorithm either makes progress toward constructing a reconfiguration $S$--$R$ path or produces a small dominating set. The $r$-claw-free condition on $(G,\mathcal{V})$ is needed to bound the number of extra vertices that contribute the error term (see Section~\ref{discussion} for further discussion about this condition).  The running time of the algorithm is of the form $O\bigl(m^{f(r,\epsilon)}n^2\bigr)$, where $n=|V(G)|$ and $f(r,\epsilon)$ depends only on $r$ and $\epsilon$. In particular, when $r$ and $\epsilon$ are fixed constants, the algorithm runs in polynomial time in $n$ and $m$. A precise running-time analysis is given in Section~\ref{sec:proof:runningtime}.

Our main theorem has consequences for the quantitative and algorithmic questions raised by Buys, Kang, and Ozeki and Wdowinski discussed above. Whenever the second outcome of Theorem~\ref{thm:main} is excluded by an appropriate domination condition, the algorithm constructs, in polynomial time, a polynomial-length reconfiguration path between any two ITs. Thus, in these settings, the reconfiguration graph $\textbf{RG}(G, \mathcal{V})$ has polynomial diameter and the corresponding reconfiguration problem can be solved efficiently. In particular, for $r$-claw-free graphs, this gives an algorithmic version of Theorem~\ref{wdowinski}(a) of Wdowinski, with  $2|I|+1$  replaced by $(2+\epsilon)|I|$. It also provides the following  answer to the question of Buys, Kang, and Ozeki (and Wdowinski) concerning the diameter, for bounded maximum degree, by taking $\epsilon<1/\Delta$ and  $r=\Delta+1$ (using the simple fact that a set $D$ can dominate at most $|D|\Delta$ vertices).

\begin{corollary} Let $\Delta$ be fixed.
Let $G$ be a graph of maximum degree $\Delta$, and let $\mathcal{V}=(V_1,\ldots,V_m)$ be a partition of $V(G)$ such that $|V_i|\geq 2\Delta+1$ for every $i\in[m]$. Then the diameter of $\textbf{RG}(G, \mathcal{V})$ is bounded by a polynomial function of $|V(G)|$.
\end{corollary}

Thus, in the bounded degree case, the same (best possible) assumption on the minimum size of the partition classes that guarantees connectivity of $\textbf{RG}(G, \mathcal{V})$ also guarantees that it has polynomial diameter.

The structural information in the second outcome of Theorem~\ref{thm:main} also allows us to obtain a concrete application to hypergraph matchings. An $r$-uniform \textit{bipartite} hypergraph $H=(A,B,E)$ on the vertex set $A\cup B$ (where $A\cap B=\emptyset$) is a hypergraph in which each edge $e\in E$ contains one vertex in $A$ and $(r-1)$ vertices in $B$, and a \textit{perfect matching} in $H$ is a set of $|A|$ disjoint edges. Given two perfect matchings $\mathcal{M}_1$ and $\mathcal{M}_2$, our algorithm finds in polynomial time either a polynomial-length reconfiguration $\mathcal{M}_1$--$\mathcal{M}_2$ path or returns a subset $S\subseteq A$ whose link $E_S$ in $G$ has a small cover. Consequently, whenever no such $S$ exists, any two perfect matchings can be connected by a polynomial-length reconfiguration path, which can also be found in polynomial time. This application is developed in Section~\ref{sec:Application}.

The remainder of the paper is organized as follows. In Section~\ref{sec:Def}, we introduce the definitions and notation used throughout the paper. In Section~\ref{sec:algorithm} we describe the algorithm underlying our main result, and establish two of its basic properties in Section~\ref{sec:proof:analysis}. In Section~\ref{sec:proof:dominatingset}, we explain how we obtain the second outcome of the main theorem. In Section~\ref{sec:proof:SRPath}, we explain how we obtain the first outcome of the main theorem and prove Theorem~\ref{thm:main}, postponing the running-time analysis to Section~\ref{sec:proof:runningtime}. In Section~\ref{sec:Application}, we present the application of our main theorem to the reconfiguration of hypergraph matchings. We end with a brief section on conclusions, discussion and open problems.

\section{Definitions and Preliminary Results}\label{sec:Def}

This section introduces the key terminology and concepts used throughout our arguments. For every nonnegative integer $n$, we extend the standard notation $[n]:=\{1,\dots, n\}$ to $[n]_0:=[n]\cup \{0\}$. For a graph $G$ with a vertex partition $\mathcal{V}=(V_1, V_2,\dots,V_m)$, a 
\textit{partial independent transversal} (PIT) is an independent set $M$ in $G$ that contains at most one vertex from each vertex class. For a set $\mathcal{B}$ of vertex classes, $M$ is an \textit{IT of $\mathcal{B}$} if $M$ is a PIT of $G$ that intersects every class in $\mathcal{B}$ and no other vertex class. For the following definitions, let $G$ be a graph with a vertex partition $\mathcal{V}=(V_1,V_2,\dots,V_m)$, and let $M$ be a PIT of $G$.

\begin{definition}\label{def:cls}
For a vertex $v\in V(G)$, let $\operatorname{cls}(v)$ denote the unique vertex class containing $v$. More generally, for any nonempty set $W\subseteq V(G)$, define $\operatorname{cls}(W):=\bigcup_{v\in W} \operatorname{cls}(v)$ and $\operatorname{ind}(W):=\{i\ : V_i=\operatorname{cls}(v)\ \emph{for} \ v\in W\}$. For an index set $I\subseteq [m]$, define $\mathcal{V}_{I}:=\{V_i:\ \emph{for} \ i\in I\}$ that is, the collection of vertex classes indexed by $I$. 
\end{definition}

\begin{definition}\label{biconstellation}
Let $\mathcal{B}$ be a set of partition classes of $(G,\mathcal{V})$.

\begin{enumerate}
\item[-] A \emph{biconstellation} for $\mathcal{B}$ is an induced subgraph $C$ of $G_{\mathcal{B}}$ such that:
\begin{enumerate}
    \item[(1)] $C$ contains a balanced bipartite graph $F$ with parts $\mathcal{P}_1$ and $\mathcal{P}_2$ (where $|\mathcal{P}_1|=|\mathcal{P}_2|$) as an induced subgraph, and every connected component of $C\setminus F$ is a star with at least two vertices.
    \item[(2)] If $L$ denotes the union of the sets of leaves of the star components, then both $L\cup \mathcal{P}_1$ and $L\cup \mathcal{P}_2$ are IT's of $\mathcal{B}$.
\end{enumerate}

\item[-] A \emph{constellation} for $\mathcal{B}$ is an induced subgraph $C$ of $G_{\mathcal{B}}$ such that:
\begin{enumerate}
    \item[(1)] Every connected component of $C$ is a star with at least two vertices.
    \item[(2)] The union of the sets of leaves of the star components forms an IT of $|\mathcal{B}|-1$ vertex classes of $\mathcal{B}$.
\end{enumerate}
\end{enumerate}
\end{definition}

\begin{figure}[t!]
    \centering   \includegraphics[scale=1.1]{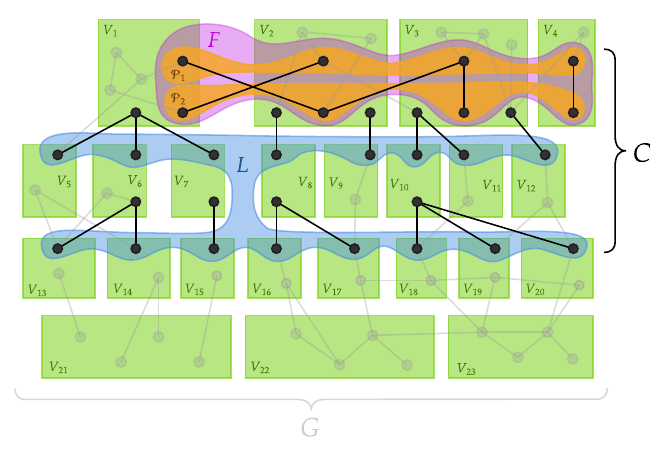}
    \caption{A graph $G$ with a vertex partition $\mathcal{V}=(V_1,V_2,\ldots,V_{23})$ and a biconstellation $C$ on the set of partition classes $\{V_1,V_2,\ldots,V_{20}\}$, where $F$ is the balanced bipartite graph of $C$ with parts $\mathcal{P}_1$ and $\mathcal{P}_2$, and $L$ is  the union of the sets of leaves of the star components.} 
    \label{fig:allsix}
\end{figure}

\begin{definition}\label{blocks}
A vertex $u$ \emph{blocks} a vertex $v$ if $u\in M$ and $uv\in E(G)$.
\end{definition}

\begin{definition}\label{immaddable}
A vertex $v$ is \emph{immediately addable} with respect to $M$ if $v\notin M$ and $N(v)\cap M=\emptyset$. For a set $W\subseteq V(G)$, we write $I_M(W)$ for the set of immediately addable vertices of $W$ with respect to $M$.
\end{definition}

\begin{definition}\label{feasible} Fix $r\ge 2$ and $\epsilon>0$. A triple $(\mu,U,\rho)$ of positive real numbers is said to be \emph{feasible} for $(r,\epsilon)$ provided that the following three inequalities hold: \begin{enumerate} 
\item[\textnormal{e\textsubscript{1}:}]\label{e_1} $\frac{2+\mu(r+2)+\rho(r+1)}{1-\mu}<(2+\epsilon)\left[1-\frac{1}{U}\left(\frac{1+\mu U}{1-\mu}+\rho\right)\right],$
\vspace{0.5em}
\item[\textnormal{e\textsubscript{2}:}]\label{e_2} $\frac{\mu(r+4)+\rho(r+2)}{1-\mu}< \epsilon\left[1-\frac{1}{U}\left(\frac{1+\mu U}{1-\mu}+\rho\right)\right],$
\vspace{0.5em}
\item[\textnormal{e\textsubscript{3}:}]\label{e_3} $\rho(1+\mu)< U.$ 

\end{enumerate} 
\end{definition}

Throughout the paper, we assume that $\mu=\rho=\frac{1}{U}=\frac{\epsilon}{10r}$. One can readily verify that $(\mu,U,\rho)=(\frac{\epsilon}{10r},\frac{10r}{\epsilon},\frac{\epsilon}{10r})$ is feasible for $(r,\epsilon)$ when $r\geq 2$ and $0<\epsilon<1$.

\begin{definition}
    A \emph{layer} with respect to the PIT $M$ is a pair $(X,Y)$ satisfying the following conditions:
    \begin{enumerate} \item $X\subseteq V(G)\setminus M$ is an independent set; \item $Y\subseteq M$ consists precisely of the vertices of M that block X; \item each vertex of Y has exactly one neighbor in X. \end{enumerate}
\end{definition}

\begin{definition}\label{tree} Consider a vertex class $A$ of $(G,\mathcal{V})$ that does not contain a vertex in $M$. Then an  \emph{alternating tree} $T$ with respect to $M$ and $A$ is a tuple $(L_0,\dots,L_\ell)$ where $\ell\ge0$ such that:
\begin{enumerate}
\item $L_0=(X_0,Y_0)=(\emptyset,\emptyset)$ and $\operatorname{cls}(Y_0):=A$,
\item $L_i=(X_i,Y_i)$ is a layer for each $i\in [\ell]$,
\item $X_i\subseteq \operatorname{cls}(Y_{i-1})$ for all $i\in [\ell]$, and
\item $(X_i\cup Y_i)\cap (X_{i'}\cup Y_{i'})=\emptyset$ for all distinct $i,i'\in[\ell]_0$.
\end{enumerate}
We call $A$ the {\it root} of $T$. 
\end{definition}

\begin{figure}[t!]
    \centering   \includegraphics[scale=1.1]{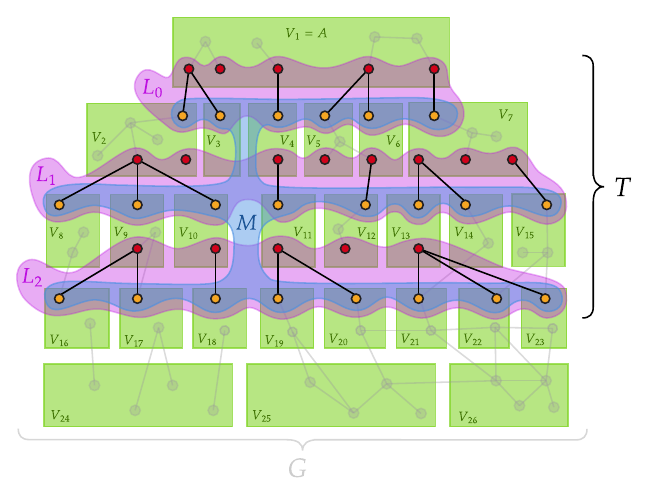}
    \caption{A graph $G$ with a vertex partition $\mathcal=(V_1,V_2,\ldots,V_{26})$ and an alternating tree $T$ rooted at $A$ with respect to a PIT $M$. The alternating tree $T$ consists of layers $L_i=(X_i,Y_i)$, for $i\in[2]_0$, where the vertices of $X_i$ and $Y_i$ are shown in red and orange, respectively.} 
    \label{fig:allsix}
\end{figure}

\begin{definition}\label{addabledef}
Let $T=(L_0,\dots,L_\ell)$ be an alternating tree of $G$ with respect to the PIT $M$ and a root $A$ and let $X,Y\subseteq V(G)$. A vertex $v\in \operatorname{cls}(Y_\ell)$ is said to be  \emph{addable} for the triple $(X,Y,T)$ if the following conditions hold; 

\begin{enumerate}
\item $v\notin Y_{\ell}\cup X\cup Y$, 

\item $N_G(v)\cap( X_{\le\ell}\cup Y_{\le\ell}\cup X\cup Y)=\emptyset$,

\item $|\operatorname{cls}(v)\cap X|<U$ whenever $\ell\ge 1$.

\end{enumerate}

\begin{definition}\label{collapsible} A layer $L_i=(X_i,Y_i)$ is \emph{collapsible} if $I_M(X_i)>\mu|X_i|$. 
\end{definition}

\begin{definition}\label{sig}
Let $T=(L_0,\dots,L_\ell)$ be an alternating tree with respect to a PIT $M$ and a vertex class $A$. For each layer $L_i$, define its \emph{signature} by $$(s_{2i-1},s_{2i})=\left(-\left\lfloor \log_b \frac{\rho^{-i}}{(1-\mu)^{i-1}}|X_i|\right\rfloor, \left\lfloor \log_b \frac{\rho^{-i}}{(1-\mu)^i}|Y_i|\right\rfloor\right)$$ where $b=\frac{U}{U-\mu\rho}$.
The {\it signature vector} of $T$ is ${s=(s_1,s_2,\dots,s_{2\ell-1}, s_{2\ell},\infty)}$.
\end{definition}

\end{definition}

\section{Algorithm}\label{sec:algorithm}
For the algorithms and the proof of Theorem~\ref{thm:main}, we switch to primed notation for the objects appearing in its statement, writing $G'$, $\mathcal V'$, $B'$, $D'$, and $C'$ in place of $G$, $\mathcal V$, $B$, $D$, and $C$, respectively. This allows us to reserve the corresponding unprimed notation for auxiliary objects introduced in the analysis. For convenience, we restate Theorem~\ref{thm:main} at the beginning of Section~\ref{sec:proof:analysis} using this primed notation.

\subsection{BuildLayer}
We begin with the basic subroutine \textsc{BuildLayer}. Let $M$ be a PIT of a pair $(G,\mathcal{V})$, and let $T=(L_0,\ldots,L_\ell)$ be an alternating tree with respect to $M$, where $L_i=(X_i,Y_i)$ for each $i\in[\ell]_0$. Given two vertex sets $X,Y\subseteq V(G)$, the subroutine \textsc{BuildLayer}$(T,X,Y)$ extends $(X,Y)$ by repeatedly adding to $X$ vertices that are addable for $(X,Y,T)$, together with their blockers in $M$ to $Y$. It terminates when no further vertex can be added in this way and returns the resulting pair $(X,Y)$.

\begin{algorithmic}[1]
\Function{BuildLayer}{$T,X,Y$}
\While{there is a vertex $v\in \operatorname{cls}(Y_\ell)$ that is addable for $(X,Y,T)$}
	\State $X:=X\cup \{v\}$
	\State $Y:=Y\cup\{u\in M:uv\in E(G)\}$
\EndWhile
\State \Return $(X,Y)$
\EndFunction
\end{algorithmic}

\subsection{SuperposedBuild}
The next subroutine, \textsc{SuperposedBuild}, is used to rebuild the layers of an existing alternating tree. Its input consists of an alternating tree $T=(L_0,\ldots,L_\ell)$, where $L_i=(X_i,Y_i)$, together with its last layer index $\ell$. Starting from the first layer, the subroutine applies \textsc{BuildLayer} to determine whether a layer can be enlarged significantly. Whenever the rebuilt version of a layer $X_i$ has size at least $(1+\mu)|X_i|$, the layer is replaced by its enlarged version and all subsequent layers are discarded. The subroutine returns the resulting alternating tree and the index of its last layer.

\begin{algorithmic}[1]
\Function{SuperposedBuild}{$T,\ell$}
\State $i:=1$
\While{$i\le\ell$}
	\State $(X'_{i},Y'_{i}):=\textsc{BuildLayer}((L_0,\dots,L_{i-1}),X_{i},Y_{i})$
	\If{$|X'_{i}|\ge(1+\mu)|X_{i}|$}
		\State $(X_{i},Y_{i}):=(X'_{i},Y'_{i})$
		\State $L_i:=(X_i,Y_i)$
		\State $T:=(L_0,\dots,L_i)$
		\State $\ell:=i$
	\EndIf
	\State $i:=i+1$
\EndWhile
\State \Return $(T,\ell)$
\EndFunction
\end{algorithmic}

\subsection{GrowPath}
We next describe the main subroutine \textsc{GrowPath}. Its input is a pair of independent transversals $\hat S$ and $\hat R$. The subroutine starts with their common part $M=\hat S\cap\hat R$ and constructs the auxiliary graph, vertex classes, and root used to build an alternating tree $T$. It then grows $T$ one layer at a time using \textsc{BuildLayer}. When sufficiently many immediately addable vertices are found, $M$ is modified and the previous layers are updated using \textsc{SuperposedBuild}; the corresponding changes to $\hat S$ and $\hat R$ are recorded as reconfiguration steps. The subroutine terminates in one of two ways. Either it makes progress towards bringing $\hat S$ and $\hat R$ together and returns the resulting reconfiguration paths $P_{\hat{S}}$ and $P_{\hat{R}}$, or it encounters a layer that is too small to continue the growth process. In the latter case, it returns the alternating tree together with a flag indicating this outcome; this tree will be used by the main algorithm to construct the sets $\mathcal{B}$ and $D$ appearing in the second outcome of Theorem~\ref{thm:main}. In this subroutine, the following notation is used.

\begin{definition}
    For a vertex $v$, let $N(v)$ denote the set of its neighbors, and let $N[v]=N(v)\cup\{v\}$. For a vertex set $S$, define $N(S)=\bigcup_{v\in S}N(v)$ and $N[S]=N(S)\cup S$.
\end{definition}

\begin{definition}
    For a vertex $v$, let $P_v$ denote a path starting at $v$, and let $\operatorname{End}(P_v)$ denote its other endpoint. Also, for a vertex $u\notin V(P_v)$ adjacent to $\operatorname{End}(P_v)$, we denote by $P_v+u$ the path obtained from $P_v$ by appending $u$ to $\operatorname{End}(P_v)$.
\end{definition}

\begin{definition}
For an IT $S$ of $G$ and a vertex $v\in V_i$ for some $i$, if $v$ has no neighbour in $S$, then \emph{$S\oplus v$} denotes the IT obtained from $S$ by replacing the vertex of $S$ in $V_i$ with $v$.
\end{definition}

\begin{algorithmic}[1]
\Function{GrowPath}{$\hat{S},\hat{R}$}
\State $M:=\hat{S}\cap \hat{R}$ 
\State $I:=\operatorname{ind}(\hat{S}\Delta \hat{R})$
\State $G:=G'\setminus N[\hat{S}\Delta \hat{R}]$
\State $V_i:=V'_i\setminus N[\hat{S}\Delta \hat{R}] $ for all $i\in[m]$
\State $A=V_0:=\cup_{i\in I} V_i$
\State $L_0:=(X_0,Y_0):=(\emptyset, \emptyset)$
\State $\operatorname{cls}(Y_0):=A$
\State $T:=(L_0)$
\State $\ell:=0$
\State $(P_{\hat{S}},P_{\hat{R}}):=(\hat{S}, \hat{R})$
\If{$G'[\hat{S}\Delta \hat{R}]$ has an isolated vertex $v$}
		\State $\Tilde{S}:= ((\hat{S}\setminus \hat{R}) \cup M)\oplus v$
        \State $\Tilde{R}:= ((\hat{R}\setminus \hat{S}) \cup M)\oplus v$ 
        \State $P_{\hat{S}}:=P_{\hat{S}} + \Tilde{S}$
        \State $P_{\hat{R}}:=P_{\hat{R}} + \Tilde{R}$ 
		\State \Return $(M,T,0,P_{\hat{S}},P_{\hat{R}})$ and terminate
\EndIf
\If{$A=\emptyset$}
    \State \Return $(M,T,2,P_{\hat{S}},P_{\hat{R}})$ and terminate
\EndIf

\While{$A\cap M=\emptyset$}
	\State $(X_{\ell+1},Y_{\ell+1}):=\textsc{BuildLayer}(T, \emptyset, \emptyset)$
	\State $L_{\ell+1}:=(X_{\ell+1},Y_{\ell+1})$
	\State $T:=(L_0,\dots,L_\ell,L_{\ell+1})$
	
	\If{$|X_{\ell+1}|\le\rho|Y_{\le\ell}|$}
		\State \Return $(M,T,1,P_{\hat{S}},P_{\hat{R}})$ and terminate
	\Else
		\State $\ell:=\ell+1$

		\While{$|I_M(X_\ell)|>\mu|X_\ell|$}
			\If{$\ell=1$}
				\State For some arbitrary $v\in I_M(X_1)$, set
                \State $\Tilde{S}:= ((\hat{S}\setminus \hat{R}) \cup M)\oplus v$
                \State $\Tilde{R}:= ((\hat{R}\setminus \hat{S}) \cup M)\oplus v$ 
                \State $P_{\hat{S}}:=P_{\hat{S}} + \Tilde{S}$
                \State $P_{\hat{R}}:=P_{\hat{R}} + \Tilde{R}$
                
				\State \Return $(M,T,0,P_{\hat{S}},P_{\hat{R}})$ and terminate
			\Else
				\ForAll{$w\in Y_{\ell-1}$ such that $I_M(X_\ell)\cap \operatorname{cls}(w)\neq\emptyset$}
					\State $M:=(M\setminus\{w\})\cup\{u\}$ for some arbitrary ${u\in I_M(X_\ell)\cap \operatorname{cls}(w)}$
					\State $Y_{\ell-1}:=Y_{\ell-1}\setminus\{w\}$
                    \State $\Tilde{S}:= ((\hat{S}\setminus \hat{R}) \cup M)$
                    \State $\Tilde{R}:= ((\hat{R}\setminus \hat{S}) \cup M)$ 
                    \State $P_{\hat{S}}:=P_{\hat{S}} + \Tilde{S}$
                    \State $P_{\hat{R}}:=P_{\hat{R}} + \Tilde{R}$

				\EndFor
			\EndIf
			\State $T':=(L_0,\dots,L_{\ell-1})$
			\State $\ell':=\ell-1$
			\State $(T,\ell):=\textsc{SuperposedBuild}(T',\ell')$			
		\EndWhile
	\EndIf
\EndWhile
\EndFunction
\end{algorithmic}

\subsection{Main Algorithm}

We are now ready to describe the main algorithm, \textsc{FindSR-Path}. Its input is an $r$-claw-free graph $G'$ with vertex partition $\mathcal V'=(V'_1,\ldots,V'_m)$ and two independent transversals $S$ and $R$. Starting from $S$ and $R$, the algorithm repeatedly calls \textsc{GrowPath} and concatenates the reconfiguration paths returned by these calls. As long as \textsc{GrowPath} makes progress, this process reduces the symmetric difference between the current pair of independent transversals. If the two transversals eventually coincide, the accumulated paths are combined to give a reconfiguration path from $S$ to $R$. Otherwise, \textsc{GrowPath} returns the exceptional outcome, from which \textsc{FindSR-Path} constructs the dominating set appearing in the second outcome of Theorem~\ref{thm:main}.
\begin{algorithmic}[1]
\Function{FindSR-Path}{$G',\mathcal{V'}=(V'_1,\dots,V'_m), S,R$}
\State $(P_S,P_R):=(\emptyset,\emptyset)$
\State $(\hat{S},\hat{R}):=(S,R)$
\While{$\hat{S}\Delta \hat{R}\neq \emptyset$}
	\State $(M,T,x,P_{\hat{S}},P_{\hat{R}}):=\textsc{GrowPath}(\hat{S},\hat{R})$
    \State $P_S:=P_S+P_{\hat{S}}$
    \State $P_R:=P_R+P_{\hat{R}}$
    \State $\hat{S}:=\text{End($P_{\hat{S}}$)}$
    \State $\hat{R}:=\text{End($P_{\hat{R}}$)}$
    \State $I:=\operatorname{ind}(\hat{S}\Delta \hat{R})$
    \If{$x=2$}
	\State $\mathcal{B}':= \mathcal{V}_I$		
	\State $D':=\hat{S}\Delta \hat{R}$
	\State \Return $\mathcal{B'}$, $D'$ and \lq\lq $G'_{\mathcal B'}$ is dominated by $D'\ $\rq\rq $\hspace{0.1em}$ and terminate
    
    \EndIf
    
    \If{$x=1$}
		\State ${\mathcal B}:=\mathcal{V}_{\operatorname{ind}(Y_{\leq\ell})}$
		\For{$j=1,\dots,\ell$}
			\State $(X'_j,Y'_j):=\textsc{BuildLayer}((L_0,\dots,L_{j-1}),X_j,Y_j)$
			\EndFor
        \State ${\mathcal{B}}:=\mathcal{ B}\setminus A^U$ as $(\hyperref[U_classes]{\text{i}})$ in  Lemma~\ref{domin_set_G}
        \State ${\mathcal B}:={\mathcal B}\setminus\mathcal{V}_{J}$ where $J=\left(\bigcup\limits_{i=2}^\ell \operatorname{ind}(X'_i\setminus X_i)\right)$ as $(\hyperref[Xi'_Xi(2_L)]{\text{ii}})$ in  Lemma~\ref{domin_set_G}
		\State $\mathcal{B}':=(\mathcal{B}\setminus \{V_0\})\cup \mathcal{V}_I$ as $(\hyperref[domin_set_B'_in_G']{\text{i}})$ in Lemma~\ref{domin_set_G'}
		
		\State $D':=D\cup (\hat{S}\Delta \hat{R})=(X'_{\le\ell}\cup Y'_{\le\ell}\cup X_{\ell+1}\cup Y_{\ell+1}\cup N)\cup (\hat{S}\Delta \hat{R})$ as $(\hyperref[domin_set_D'_in_G']{\text{ii}})$ in Lemma~\ref{domin_set_G'}
		\State \Return $\mathcal{B'}$, $D'$ and \lq\lq $G'_{\mathcal B'}$ is dominated by $D'\ $\rq\rq $\hspace{0.1em}$ and terminate
	\EndIf
\EndWhile
\State \Return $P=P_S+P_R$ as  \lq\lq a $S$--$R$ path in $\textbf{RG}(G,\mathcal{V})$\rq\rq
\EndFunction
\end{algorithmic}\vspace{6pt}

\section{Proofs: Analysis}\label{sec:proof:analysis}

Following the explanation at the beginning of Section~\ref{sec:algorithm}, we restate Theorem~\ref{thm:main} below using the primed notation.
\begin{mainthmprime}\label{thm:mainprime}
Let $r\in \poi$ and $\epsilon >0$ be given. Then there is an algorithm $A_{r,\epsilon}$ that takes an $r$-claw-free graph $G'$ with respect to a vertex partition $\mathcal{V'}=(V'_1, V'_2,\dots,V'_m)$ and two independent transversals $S,R\in \mathcal{I}(G',\mathcal{V'})$ as input, and finds in polynomial time either;
    \begin{itemize}
        \item [-] a polynomial-length $S$--$R$ path $P$ in $\textbf{RG}(G', \mathcal{V'})$, or
            \item[-] a set of partition classes $\mathcal{B'}\subseteq \mathcal{V}$ and a vertex set $D'\subseteq V(G')$ such that $|D'|\leq (2+\epsilon)|\mathcal{B'}|$ and $D'$ dominates $G'_{\mathcal{B'}}$. Furthermore, there exists a biconstellation $C'$ for some $\mathcal{B}'_0\supseteq \mathcal{B'}$ where with $V(C')\subseteq D'$ and $|D'\setminus V(C')|< \epsilon(|\mathcal{B'}|-1)$.
    \end{itemize}
\end{mainthmprime}

First, as a preprocessing step, we establish via the following easy lemma that we may assume the input graph in Theorem~$\hyperref[thm:mainprime]{1.4}$ has no isolated vertices.

\begin{lemma}\label{lem:isolatedvertex} 
For $r\in\poi$ and $\epsilon>0$, let $G',\mathcal{V}',S,$ and $R$ be as in \emph{Theorem~$\hyperref[thm:mainprime]{1.4}$}, and let $v\in V'_i$ be an isolated vertex for some $i\in[m]$. Let $v_S$ and $v_R$ be the unique vertices of $S\cap V'_i$ and $R\cap V'_i$, respectively. If $A_{r,\epsilon}$ is an algorithm satisfying \emph{Theorem~$\hyperref[thm:mainprime]{1.4}$} for the input $(G'-V'_i,\mathcal{V}'\setminus\{V'_i\},S\setminus\{v_S\},R\setminus\{v_R\})$, then there exists an algorithm $\hat{A}_{r,\epsilon}$ satisfying \emph{Theorem~$\hyperref[thm:mainprime]{1.4}$} for the input $(G',\mathcal{V}',S,R)$. 
\end{lemma} 
\begin{proof} 
We define $\hat{A}_{r,\epsilon}$ as follows. Given the input $(G',\mathcal{V}',S,R)$, the algorithm $\hat{A}_{r,\epsilon}$ runs $A_{r,\epsilon}$ on $(G'-V'_i,\mathcal{V}'\setminus\{V'_i\},S\setminus\{v_S\},R\setminus\{v_R\})$ and extends its output to an output for $(G',\mathcal{V}',S,R)$, as described below. Suppose first that $A_{r,\epsilon}$ returns a set of partition classes $\mathcal{B}'$ and a vertex set $D'$ satisfying the second outcome of Theorem~$\hyperref[thm:mainprime]{1.4}$. Since $\mathcal{B}'\subseteq \mathcal{V}'\setminus\{V'_i\}\subseteq\mathcal{V}'$ and $D'\subseteq V(G'-V'_i)\subseteq V(G')$, both $\mathcal{B}'$ and $D'$ satisfy the required containment conditions for $G'$. Moreover, since $V'_i\notin \mathcal{B}'$, we have $(G'-V'_i)_{\mathcal{B}'}=G'_{\mathcal{B}'}$, and hence $D'$ dominates $G'_{\mathcal{B}'}$. All other conditions, including the bound on $|D'|$ and those concerning the biconstellation, remain unchanged. Thus, $\mathcal{B}'$ and $D'$ satisfy the second outcome of Theorem~$\hyperref[thm:mainprime]{1.4}$ for $(G',\mathcal{V}',S,R)$. Now, suppose that $A_{r,\epsilon}$ returns a polynomial-length path $P=(S\setminus\{v_S\}=I_0,I_1,\ldots,I_t=R\setminus\{v_R\})$ in $\textbf{RG}(G'-V'_i,\mathcal{V}'\setminus\{V'_i\})$. The algorithm $\hat{A}_{r,\epsilon}$ adds $v$ to every node of $P$. Since $v$ is isolated, each $I_j\cup\{v\}$ is an independent transversal for $(G',\mathcal{V}')$. Moreover, these sets form a path in $\textbf{RG}(G',\mathcal{V}')$, since consecutive sets differ in exactly one vertex. If its starting node is not $S$, we prepend $S$ to the path by replacing $v$ with $v_S$; similarly, if its ending node is not $R$, we append $R$ by replacing $v$ with $v_R$. The resulting path is an $S$--$R$ path in $\textbf{RG}(G',\mathcal{V}')$, and its length is at most two more than the length of $P$, so it is still polynomial. Finally, since $A_{r,\epsilon}$ runs in polynomial time and the additional steps above clearly take polynomial time, $\hat{A}_{r,\epsilon}$ also runs in polynomial time. Therefore, $\hat{A}_{r,\epsilon}$ satisfies Theorem~$\hyperref[thm:mainprime]{1.4}$ for the input $(G',\mathcal{V}',S,R)$. 
\end{proof}

Hence, by Lemma~\ref{lem:isolatedvertex}, we may assume throughout the proof of Theorem~$\hyperref[thm:mainprime]{1.4}$ that the input graph has no isolated vertices. Now, to prove Theorem~$\hyperref[thm:mainprime]{1.4}$, we divide the proof into four sections. In this section, we state two properties that will be used in the proof of Lemma~\ref{domin_set_G} in the next section. In Lemma~\ref{no_collapse}, we establish a basic property of the layers maintained by \textsc{GrowPath}, and in Lemma~\ref{nosuper}, we show that no layer produced by \textsc{BuildLayer} can be substantially enlarged.
Hence, by Lemma~\ref{lem:isolatedvertex}, we assume during the proof of  Theorem~$\hyperref[thm:mainprime]{1.4}$, we divide the proof into four sections. In this section, we state two properties that will be used in the proof of Lemma~\ref{domin_set_G} in the next section. In Lemma~\ref{no_collapse}, we establish a basic property of the layers maintained by \textsc{GrowPath}, and in Lemma~\ref{nosuper}, we show that no layer produced by \textsc{BuildLayer} can be substantially enlarged.

\begin{lemma}\label{no_collapse} Let $T=(L_0,\dots,L_\ell)$ and $M$ be the alternating tree and partial independent transversal at the beginning of an iteration of the while loop in line $22$ of {\rm \textsc{GrowPath}}. Then none of the layers $L_0,\dots,L_\ell$ are collapsible. Consequently, $|Y_i|\ge (1-\mu)|X_i|$ for all $i\in[\ell]$. 
\end{lemma} 

\begin{proof} We prove the first statement by induction on the iterations of the while loop in line~22. The claim is true at the beginning of the first iteration, since the alternating tree initially consists only of the root layer where $|I_M(X_i)|=\mu|X_i|=0$. \\
Now assume that, at the beginning of the current iteration, none of the layers of the alternating tree are collapsible. During lines~$23$--$24$, the algorithm constructs a new layer $L_{\ell+1}$. If {\rm \textsc{GrowPath}} terminates during this iteration, there is nothing to prove. Otherwise, there are two cases:

\textbf{Case 1:} $L_{\ell+1}$ is not collapsible. Then no collapse operation is performed. The previously existing layers are not modified, so they remain non-collapsible. Since the new layer is also non-collapsible, the claim holds at the beginning of the next iteration.

\textbf{Case 2:} $L_{\ell+1}$ is collapsible. Let $T=(L_0,\dots,L_k),$ where $k\le\ell$, denote the alternating tree after the loop of collapsing operations (lines~30--51) has finished. Since the algorithm does not terminate, we have $k\ge1$. Each collapse operation (lines~39--46) removes the current last layer (line~48), thereby strictly decreasing the number of layers. Moreover, if {\rm \textsc{SuperposedBuild}} modifies the tree in line~50, it either preserves or decreases the number of layers. Consequently, the layers $L_0,\dots,L_{k-1}$ are never modified during the collapsing loop and therefore remain non-collapsible by the induction hypothesis. The only layer that may change is the final layer $L_k$, through lines~40, 41, and~50. However, when the collapsing loop terminates, its stopping condition guarantees that the last layer is no longer collapsible. Hence every layer of the resulting alternating tree is non-collapsible, establishing the induction step. Therefore, at the beginning of every iteration of the while loop, none of the layers of the alternating tree are collapsible.

Finally, fix a layer $L_i$ with $ i\in[\ell]$. Since $L_i$ is not collapsible, $|I_M(X_i)|\le\mu|X_i|.$ Moreover, by the definition of a layer, every vertex of $X_i\setminus I_M(X_i)$ has a distinct blocking vertex in $Y_i$. Hence $|X_i\setminus I_M(X_i)|\le |Y_i|$. Combining these inequalities yields \[ |X_i| =|X_i\setminus I_M(X_i)|+|I_M(X_i)| \le |Y_i|+\mu|X_i|, \] or equivalently, $ |Y_i|\ge(1-\mu)|X_i|, $ as claimed. 
\end{proof}

\begin{lemma}\label{nosuper} Let $T=(L_0,\dots,L_\ell)$ and $M$ be the alternating tree and partial independent transversal at the beginning of an iteration of the while loop in line~22 of {\rm \textsc{GrowPath}}. For every $i\in[\ell]$, define $(X'_i,Y'_i):=\emph{\textsc{BuildLayer}}((L_0,\dots,L_{i-1}),X_i,Y_i)$. Then $ |X'_i|<(1+\mu)|X_i|$ for all $i\in[\ell]$.
\end{lemma} 

\begin{proof} Fix $i\in[\ell]$, and consider the layer $L_i$ at the beginning of the current iteration. When $L_i$ was created, either by {\rm \textsc{BuildLayer}} (line~23) or by {\rm \textsc{SuperposedBuild}} (line~50), the set $X_i$ was chosen to be maximal with respect to the current PIT $M$. Hence, immediately after its construction, another application of {\rm \textsc{SuperposedBuild}} could not enlarge $X_i$. If no collapsible layer has appeared since $L_i$ was created, then the PIT $M$ has remained unchanged throughout the intervening iterations. Consequently, the same maximality argument still applies, and {\rm \textsc{BuildLayer}} cannot enlarge $X_i$ by a factor of $\mu$. Now suppose that one or more collapsible layers have appeared since $L_i$ was created. Every such layer must have index larger than $i$; otherwise, the collapse operation would have removed $L_i$ from the alternating tree (lines~48--50). Whenever such a collapse occurs, {\rm \textsc{SuperposedBuild}} attempts to augment $X_i$ by at least $\mu|X_i|$. Since $L_i$ is still present in the alternating tree at the beginning of the current iteration, every one of these attempts must have failed. Therefore, even if {\rm \textsc{BuildLayer}} is executed again on the current tree, \[ (X'_i,Y'_i) = \text{\rm \textsc{BuildLayer}}((L_0,\dots,L_{i-1}),X_i,Y_i), \] it cannot produce an augmentation of size at least $\mu|X_i|$. So, $|X'_i|<(1+\mu)|X_i|,$ as required. \end{proof}

\section{Proofs: Dominating Set}\label{sec:proof:dominatingset}
In this section, we show how the second outcome of Theorem~$\hyperref[thm:mainprime]{1.4}$ arises and construct the desired dominating set and biconstellation. Assume that for two independent transversals $\hat S$ and $\hat R$, we are at the beginning of {\rm \textsc{GrowPath}($\hat S, \hat R$)}. Let $M:=\hat{S}\cap \hat{R}$, $I:=\operatorname{ind}(\hat{S}\Delta \hat{R})$, $G:=G'\setminus N[\hat{S}\Delta \hat{R}]$, $V_i:=V'_i\setminus N[\hat{S}\Delta \hat{R}] $ for all $i\in[m]$, and $A=V_0:=\cup_{i\in I} V_i$ as defined in lines $2$--$6$ of {\rm \textsc{GrowPath}}. The condition that $(G,\mathcal{V})$ is $r$-claw-free is needed in the proof of the following lemma (specifically in Claim~\ref{userclawfree}).

\begin{lemma}\phantomsection\label{domin_set_G}
Let $(\mu,U,\rho)$ be feasible for given $\epsilon >0$ and $r\in \poi$. Assume $T=(L_0,\dots,L_{\ell})$ is the alternating tree with respect to the root $A$ and PIT $M$ at the beginning of some iteration of the while loop in line $22$ of {\rm \textsc{GrowPath}}. Suppose  $|X_{\ell+1}|\leq \rho|Y_{\leq \ell}|$ in line $26$ of {\rm \textsc{GrowPath}} and the algorithm terminates in line $27$. Then there exists a set of partition classes $\mathcal{B}\subseteq (\mathcal{V}\setminus \mathcal{V}_I)\cup 
\{V_0\}$ and a vertex set $D\subseteq V(G')$ such that $|D|\leq (2+\epsilon)(|\mathcal{B}|-1)$, and $D$ dominates $G_{\mathcal{B}}$. Moreover, $D$ contains a constellation $C$ for a set of partition classes $\mathcal{B}_0\supseteq\mathcal{B}$ where $|D\setminus V(C)| < \epsilon (|\mathcal{B}|-1)$.
\end{lemma}

\begin{proof}
Let $|X_{\ell+1}|\leq \rho|Y_{\leq \ell}|$ at line $26$ of some iteration of the while loop in line $22$.
For each $i\in [\ell]$, let $(X'_i,Y'_i)=\textsc{BuildLayer}((L_0,\dots,L_{i-1}),X_i,Y_i)$. Now, obtain the set of partition classes $\mathcal{B}$ from $\mathcal{B}_0:= \mathcal{V}_{\operatorname{ind}(Y_{\leq \ell})}$ by;
\begin{enumerate}
    \item[(i)]\phantomsection\label{U_classes} removing all classes of $\mathcal{B}$ that are not $A$ and contain $U$ vertices of $X_{\leq\ell}$ (denote the set of these classes by $A^U)$, and
    \item[(ii)]\phantomsection\label{Xi'_Xi(2_L)} removing all classes in $\mathcal{V}_J$ with $J=\bigcup\limits_{i=2}^\ell \operatorname{ind}(X'_i\setminus X_i)$ from $\mathcal{B}$ (note that $\operatorname{ind}(X'_1\setminus X_1)=\emptyset$). 
\end{enumerate}

Let $W=X'_{\le\ell}\cup Y'_{\le\ell}\cup X_{\ell+1}\cup Y_{\ell+1}$ and $N\subseteq V(G')$ be the smallest set containing at least one neighbor for each vertex of $I_M(W)$ in $G'$. Such a neighbor exists for every vertex since, by Lemma~\ref{lem:isolatedvertex}, we may assume that $G'$ has no isolated vertices. Define $D=W\cup N$. Now, we claim that $|D|\leq (2+\epsilon)(|\mathcal{B}|-1)$ and $D$ dominates $G_{\mathcal{B}}$. We prove our claim step by step.

\begin{claim}
Let $B=V_{\mathcal{B}}$. Then $W$ dominates $G[B\setminus I_M(W)]$.
\end{claim}

\begin{proof}\renewcommand{\qedsymbol}{$\square$} Consider a vertex $u\in B\setminus I_M(W)$. Then we have the following cases;
\begin{itemize}
    \item[-] $u\in (Y'_{\leq\ell}\cup Y_{\ell+1}) \setminus I_M(W)$. Then, according to the lines 3-4 of \textsc{BuildLayer}, $u$ is the neighbor of a vertex  $v\in (X'_{\leq\ell}\cup X_{\ell+1})$. Hence, $v$ is the neighbor of $u$ in $W$.
    \item[-] $u\in (X'_{\leq\ell}\cup X_{\ell+1}) \setminus I_M(W)$. Then $u$ has a neighbor $v\in M$; otherwise, $u\in I_M(W)$. Hence, $v\in (Y'_{\leq\ell}\cup Y_{\ell+1})$, and it is the neighbor of $u$ in $W$.
    \item[-] $u\in B \setminus W$. 
    Since ${\operatorname{cls}(u)}\in\mathcal{B}_0=\mathcal{V}_{\operatorname{ind}(Y_{\leq\ell})}$, there is an $i\in[\ell+1]$ such that $\operatorname{cls}(u)\cap Y_{i-1}\neq \emptyset$. By (i), $\operatorname{cls}(u)$ is equal to $A$ or has an intersection less than $U$ with $X_i$. Also, $u\notin (X_{\leq i}\cup Y_{\leq i}) \subseteq W$. So, if $u$ does not have a neighbor in $X_{\leq i}\cup Y_{\leq i}$, then $u$ is an addable vertex for $X_i,Y_i,$ and $(L_0,\dots,L_{i-1})$ which means; the function  $\textsc{BuildLayer}((L_0,\dots,L_{i-1}),X_i,Y_i)$ added $u$ to $X_i$ (note that $X'_i=X_i $ for all $V_i\in \mathcal{B}$), but we know $u\notin W\supseteq X_i$. Therefore, $u$ has a neighbor $v\in X_{\leq i}\cup Y_{\leq i}$. Hence, $v$ is the neighbor of $u$ in $W$.
    \end{itemize}
Thus, $G[B\setminus I_M(W)]$ is dominated by $W$.
\end{proof}
Since $N$ clearly dominates $I_M(W)$, it follows that $D=W\cup N$ dominates $G_{\mathcal{B}}$. To show that $|D|\leq (2+\epsilon)(|\mathcal{B}|-1)$, we need the following relations.

\begin{claim}\label{userclawfree}
We have the following relations: 
\begin{itemize}
    \item[\textnormal{R\textsubscript{1}:}]\label{R_1} $|\mathcal{B}|-1\geq\left[1-\frac{1}{U}\left(\frac{1+\mu
      U}{1-\mu}+\rho\right)\right]|Y_{\le\ell}|$,\vspace{0.5em} 
    \item[\textnormal{R\textsubscript{2}:}]\label{R_2} $|N|\leq |I_M(W)|<\frac{2\mu+\rho}{1-\mu}|Y_{\le\ell}|$,\vspace{0.9em} 
    \item[\textnormal{R\textsubscript{3}:}]\label{R_3} $|Y'_{\le\ell}\setminus Y_{\le\ell}|<\mu(r-1)|X_{\le\ell}|$,\vspace{0.7em}
    \item[\textnormal{R\textsubscript{4}:}]\label{R_4} $|D\setminus (X_{\le\ell}\cup Y_{\le\ell})|< \frac{\mu(r+2)+\rho(r+1)}{1-\mu}|Y_{\le\ell}|.$
\end{itemize}
\end{claim}
\begin{proof}\renewcommand{\qedsymbol}{$\square$}
    Define $A^U$ as the set of vertex classes in $\mathcal{V}_{\operatorname{ind}(Y_{\le\ell})}$ that contain $U$ vertices in
$X_{\le\ell+1}$. Then we can list the following facts:
    \begin{itemize}
        \item [f$_1$:]\label{f_1} $|A^U|\leq U^{-1}|X_{\le\ell+1}|$,\vspace{0.5em} 
        \item [f$_2$:]\label{f_2} $|\operatorname{ind}(Y_{\le\ell})|=|Y_{\le\ell}|+1$ (as  $Y_0=\emptyset$ and $Y_i$'s are disjoint),\vspace{0.5em}
        \item[f$_3$:]\label{f_3}$|X_{\le\ell}|\le\frac{1}{1-\mu}|Y_{\le\ell}|$ (by ~\cref{no_collapse}), \vspace{0.5em}
        \item[f$_4$:]\label{f_4}$|I_M(X_{\leq\ell})|\le\mu|X_{\le\ell}|$ (all layers $\leq \ell$ are not collapsible by line 30 of \textsc{GrowPath}), \vspace{0.5em}
        \item[f$_5$:]\label{f_5} $|X_{\ell+1}|\le\rho|Y_{\le\ell}|$ (the assumption of Lemma~\ref{domin_set_G}),\vspace{0.5em}
        \item[f$_6$:]\label{f_6}$I_M(W)\subseteq
        X'_{\le\ell}\cup X_{\ell+1}$ (by Definition~\ref{immaddable}),\vspace{0.5em}
        \item[f$_7$:]\label{f_7}$|X'_{\le\ell}\setminus
X_{\le\ell}|<\mu|X_{\le\ell}|$ (
by Lemma~\ref{nosuper}, we have $|X'_{\le\ell}|<(1+\mu)|X_{\le\ell}|$),
        
        \item [f$_8$:]\label{f_8} $\sum\limits_{i=1}^{\ell} |\operatorname{ind}(X'_i\setminus X_i)|\le\mu|X_{\le\ell}|$ (Using the previous fact, since $|\operatorname{ind}(X'_i\setminus X_i)|\leq|X'_i\setminus X_i|$ ).\vspace{0.5em}
    \end{itemize}
Now, according to the above facts, we can obtain the following relations:

\begin{align*}
|{\mathcal B}|&\ge|\operatorname{ind}(Y_{\le\ell})|-|A^U|-\left|\bigcup\limits_{i=1}^\ell \operatorname{ind}(X'_i\setminus X_i)\right|&\hspace{2cm}\\
&\ge(|Y_{\le\ell}|+1)-\frac{1}{U}|X_{\le\ell+1}|-\mu|X_{\le\ell}|&(\hyperref[f_2]{\text{f$_2$}},\hyperref[f_2]{\text{f$_1$}},\hyperref[f_8]{\text{f$_8$}})\\
&=(|Y_{\le\ell}|+1)-\left(\frac{1}{U}|X_{\le\ell}|+\frac{1}{U}|X_{\ell+1}|+\mu|X_{\le\ell}|\right)&\\
&\ge |Y_{\le\ell}|+1-\left[\left(\frac{1}{U}+\mu\right)|X_{\le\ell}|+\frac{\rho}{U}|Y_{\le\ell}|\right]&(\hyperref[f_5]{\text{f$_5$}})\\
&\ge |Y_{\le\ell}|+1-\left[\left(\frac{1}{U}+\mu\right)\left(\frac{1}{1-\mu}\right)|Y_{\le\ell}|+\frac{\rho}{U}|Y_{\le\ell}|\right]&(\hyperref[f_3]{\text{f$_3$}})\\
&= 1+\left[1-\frac{1}{U}\left(\frac{1+\mu U}{1-\mu}+\rho\right)\right]|Y_{\le\ell}|&
\end{align*}

and,

\begin{align*}
|I_M(W)|&=|I_M(X'_{\le\ell}\cup X_{\ell+1})|&\hspace{2.3cm}(\hyperref[f_6]{\text{f$_6$}})\\
&=|I_M(X'_{\le\ell})|+|I_M(X_{\ell+1})|&\\
&\le|I_M(X_{\le\ell})|+|I_M(X'_{\le\ell}\setminus X_{\le\ell})|+|X_{\ell+1}|&\\
&<\mu|X_{\le\ell}|+\mu|X_{\le\ell}|+|X_{\ell+1}|&(\hyperref[f_4]{\text{f$_4$}},\hyperref[f_7]{\text{f$_7$}})\\
&=2\mu|X_{\le\ell}|+|X_{\ell+1}|&\\&=\frac{2\mu}{1-\mu}|Y_{\le\ell}|+\rho|Y_{\le\ell}|.&(\hyperref[f_3]{\text{f$_3$}},\hyperref[f_5]{\text{f$_5$}})
\end{align*}

Therefore, by definition of $N$, we can obtain $|N|\leq |I_M(W)|<\frac{2\mu+\rho}{1-\mu}|Y_{\le\ell}|$. 

For the third relation, consider a vertex $v\in Y'_{\le\ell}\setminus Y_{\le\ell}$. Since, by the construction of $(X'_i,Y'_i)$, every vertex in $Y'_{\leq\ell}$ has a unique neighbor in $X'_{\leq\ell}$, there is a unique vertex $u\in X'_i$ for some $i\in[\ell]$ where $uv\in E(G)$ and $v\in  Y'_i\setminus Y_i$. If $u\in X_i$, then by line $4$ of \textsc{BuildLayer}, we should have $v\in Y_i\subseteq Y_{\leq\ell}$ which is a contradiction. Therefore, $v\in X'_{\le\ell}\setminus X_{\le\ell}$. Also, note that $Y'_{\le\ell}$ is an independent set and $G$ is a $r$-claw-free graph, meaning each vertex of $X'_{i}$ has at most $r-1$ independent neighbors in different vertex classes. Hence, $|Y'_{\le\ell}\setminus Y_{\le\ell}|\le (r-1)|X'_{\le\ell}\setminus X_{\le\ell}|$ and so by $\hyperref[f_7]{\text{f$_7$}}$, we have $|Y'_{\le\ell}\setminus Y_{\le\ell}|< \mu(r-1)|X_{\le\ell}|,$ as desired. \\

To prove the fourth relation, recall that $D=W\cup N$. Hence, we have;
\begin{align*}
\hspace{1cm}&|D\setminus
(X_{\le\ell}\cup
Y_{\le\ell})|\leq\\&|(X'_{\le\ell}\setminus
X_{\le\ell})\cup(
Y'_{\le\ell}\setminus Y_{\le\ell})\cup X_{\ell+1}\cup Y_{\ell+1}\cup
N|\leq &\\
&|X'_{\le\ell}\setminus
X_{\le\ell}|+
|Y'_{\le\ell}\setminus Y_{\le\ell}|+ |X_{\ell+1}|+|Y_{\ell+1}|+
|N|\leq &\\
&\mu|X_{\le\ell}|+\mu(r-1)|X_{\le\ell}|+\rho|Y_{\le\ell}|+\rho(r-1)|Y_{\le\ell}|+\frac{2\mu+\rho}{1-\mu}|Y_{\le\ell}|=\hspace{0.8cm}(\hyperref[f_7]{\text{f$_7$}},\hyperref[R_3]{\text{R$_3$}},\hyperref[f_5]{\text{f$_5$}},\hyperref[f_5]{\text{f$_5$}},\hyperref[R_2]{\text{R$_2$}})&\\
&\mu r|X_{\le\ell}|+\rho r|Y_{\le\ell}|+\frac{2\mu+\rho}{1-\mu}|Y_{\le\ell}|<\mu r|X_{\le\ell}|+\frac{2\mu+\rho(r+1)}{1-\mu}|Y_{\le\ell}|<\\
&\hspace{7.2cm}\frac{\mu(r+2)+\rho(r+1)}{1-\mu}|Y_{\le\ell}|,\hspace{3.151cm}(\hyperref[f_3]{\text{f$_3$}})
\end{align*}
where in the third inequality, the bound $|Y_{\ell+1}|\leq \rho(r-1)|Y_{\leq \ell}|$ follows from combining $\hyperref[f_5]{\text{f$_5$}}$ with the fact that $|Y_{\ell+1}|\leq (r-1)|X_{\ell+1}|$, since each vertex of $X_{\ell+1}$ has at most $r-1$ independent neighbours in distinct vertex classes, as $G$ is $r$-claw-free.

\end{proof}
Now, we can prove $|D|\leq (2+\epsilon)(|\mathcal{B}|-1)$ as follows;
\begin{align*}
|D|&=|
X_{\le\ell}\cup
Y_{\le\ell}|+|D\setminus
(X_{\le\ell}\cup
Y_{\le\ell})|&\\
&<\left(\frac1{1-\mu}+1\right)|Y_{\le\ell}|+\left(\frac{\mu(r+2)+\rho(r+1)}{1-\mu}\right)|Y_{\le\ell}|&\hspace{2.3cm}(\hyperref[f_3]{\text{f$_3$}},\hyperref[R_4]{\text{R$_4$}})\\&
< \left(\frac{2+\mu(r+2)+\rho(r+1)}{1-\mu}\right)|Y_{\le\ell}|
&\\&
<(2+\epsilon)\left[1-\frac{1}{U}\left(\frac{1+\mu U}{1-\mu}+\rho\right) \right]|Y_{\le\ell}|&(\hyperref[e_1]{\text{e$_1$}})\\
&< (2+\epsilon)(|\mathcal{B}|-1).&(\hyperref[R_1]{\text{R$_1$}})\\
\end{align*}

Also, if we define $C=G[X_{\leq\ell}\cup
  Y_{\leq\ell}\setminus I_M(X_{\leq\ell})]$ then clearly, according to Definition~\ref{biconstellation} and the construction of $X_{\leq\ell}$ and $
  Y_{\leq\ell}$ by \textsc{BuildLayer}, $C$ is a constellation of $\mathcal{B}_0$ in $G$, and since $I_M(X_{\le\ell})\subseteq I_M(W)$, we have:

\begin{align*}
|D\setminus V(C)|&=|(D\setminus
(X_{\le\ell}\cup
Y_{\le\ell}))\cup I_M(X_{\le\ell})|&\\&< \frac{\mu(r+2)+\rho(r+1)}{1-\mu}|Y_{\le\ell}|+\frac{2\mu+\rho}{1-\mu}|Y_{\le\ell}|&\hspace{2.7cm}(\hyperref[R_4]{\text{R$_4$}},\hyperref[R_2]{\text{R$_2$}})\\
&=\frac{\mu(r+4)+\rho(r+2)}{1-\mu}|Y_{\le\ell}|&\\&<\epsilon\left[1-\frac{1}{U}\left(\frac{1+\mu U}{1-\mu}+\rho\right)\right]|Y_{\le\ell}|&(\hyperref[e_2]{\text{e$_2$}})\\& < \epsilon(|\mathcal{B}|-1).&(\hyperref[R_1]{\text{R$_1$}})\\
\end{align*}

\end{proof}

Now, we show that, by Lemma~\ref{domin_set_G}, if at some point during the execution of \textsc{GrowPath}, the condition in line $26$ is satisfied, then the second outcome of Theorem~$\hyperref[thm:mainprime]{1.4}$ arises.

\begin{lemma}\label{domin_set_G'}
Let $(\mu,U,\rho)$ be feasible for given $\epsilon >0$ and $r\in \poi$. Assume $T=(L_0,\dots,L_{\ell})$ is the alternating tree with respect to the root $A$ and PIT $M$ at the beginning of some iteration of the while loop in line $22$ of {\rm \textsc{GrowPath}}. If $|X_{\ell+1}|\leq \rho|Y_{\leq \ell}|$ in line $26$ of {\rm \textsc{GrowPath}} and the algorithm terminates in line $27$, then the second case of \emph{\Cref{thm:main}} holds: there exists a set of partition classes $\mathcal{B'}\subseteq \mathcal{V}$ and a vertex set $D'\subseteq V(G')$ such that $|D'|\leq (2+\epsilon)|\mathcal{B'}|$, and $D'$ dominates $G'_{\mathcal{B'}}$. Moreover, $D'$ contains a biconstellation $C'$ for a set of partition classes $\mathcal{B}'_0\supseteq\mathcal{B'}$ where $|D'\setminus V(C')| < \epsilon (|\mathcal{B'}|-1)$.
\end{lemma}

\begin{proof} Consider an execution of {\rm \textsc{GrowPath}} on input $(\hat{S}$,$\hat{R})$. Let $I,G,$ and $ (V_1,\dots,V_m)$ be defined as in lines $3$--$5$. First, note that since {\rm \textsc{GrowPath}} did not terminate at either line 17 or line 20, we may assume the following:

\begin{itemize}
    \item[a$_1$:]\label{a_1} Every vertex of the induced balanced bipartite graph $G'[\hat{S}\Delta \hat{R}]$, with partition $(\hat{S}\setminus\hat{R},\hat{R}\setminus\hat{S})$, has degree at least one: Note that, without loss of generality, if there exists a vertex $v\in \hat{S}\setminus\hat{R}$ such that $N[v]\cap \hat{R}=\emptyset$, then we can define $\Tilde{S}, \Tilde{R}, P_{\hat{S}},$ and $P_{\hat{R}}$ as in lines $13$--$16$ (note that $\hat{S}=\Tilde{S}$), return $(M,T,0,P_{\hat{S}},P_{\hat{R}})$, and terminate.
    \item[ a$_2$:]\label{a_2} $A\neq\emptyset$: Note that if $A=\cup_{i\in I} V_i=\emptyset$, then we have $\cup_{i\in I} V_i\subseteq N[\hat{S}\Delta \hat{R}]$. Combining this with $(\hyperref[a_1]{\text{a$_1$}})$, we can conclude that $D':=\hat{S}\Delta \hat{R}$ is a dominating set of $G'_I$ with $|D'|=2|I|$. Moreover, it is immediate that the induced subgraph $G'[N[\hat{S}\Delta \hat{R}]]$ is the desired biconstellation for $I$ inside $D'$ (see lines 11-15 of \textsc{FindSR-Path}).
\end{itemize}

Now, let $|X_{\ell+1}|\leq \rho|Y_{\leq \ell}|$ at the beginning of some iteration of the while loop in line $22$. By $(\hyperref[a_2]{\text{a$_2$}})$, we can assume that $A\neq \emptyset$. Therefore, according to Lemma~\ref{domin_set_G}, there exists an index set $\mathcal{B}\subseteq (\mathcal{V}\setminus \mathcal{V}_I)\cup 
\{V_0\}$ and a vertex set $D\subseteq V(G')$ such that $|D|\leq (2+\epsilon)(|\mathcal{B}|-1)$, and $D$ dominates $G_{\mathcal{B}}$. Moreover, $D$ contains a constellation $C$ for $\mathcal{B}$ where $|D\setminus V(C)| < \epsilon (|\mathcal{B}|-1)$. Recall from Lemma~\ref{domin_set_G} that $D=W\cup N=(X'_{\le\ell}\cup Y'_{\le\ell}\cup X_{\ell+1}\cup Y_{\ell+1})\cup N$ where $N\subseteq V(G')$ is a smallest set containing at least one neighbor of each vertex in $I_M(W)$.

Now, consider $\mathcal{B}', D'$ and $C'$ as follows:
\begin{itemize}\phantomsection\label{domin_set_B'D'C'_in_G'}
    \item[(i)]\phantomsection\label{domin_set_B'_in_G'} $\mathcal{B}'=(\mathcal{B}\setminus \{V_0\})\cup \mathcal{V}_I$, $\mathcal{B}'_0=(\mathcal{B}_0\setminus \{V_0\})\cup \mathcal{V}_I$: we have $\mathcal{B}'\subseteq \mathcal{V}$ and $|\mathcal{B}'|=|\mathcal{B}|+|I|-1$.
    \item[(ii)]\phantomsection\label{domin_set_D'_in_G'} $D'=D\cup (\hat{S}\Delta \hat{R})=((X'_{\le\ell}\cup Y'_{\le\ell}\cup X_{\ell+1}\cup Y_{\ell+1})\cup N)\cup (\hat{S}\Delta \hat{R}): $ note that $|D'|\leq (2+\epsilon)(|\mathcal{B}|-1) +2|I|\leq (2+\epsilon)|\mathcal{B}'|$. Note that $V_0\in \mathcal{B}$ meaning $\cup_{i\in I} V_i\subseteq V(G_{\mathcal{B}})$. Combining this with $(\hyperref[a_1]{\text{a$_1$}})$, we conclude that that $D'$ dominates $G'_{\mathcal{B'}}$.
    \item[(iii)] $C'= C\cup G'[\hat{S}\Delta \hat{R}]$: let $L$ denote the union of the sets of leaves of the star components of constellation $C$. Since $L\subseteq   Y_{\leq\ell}\subseteq M=\hat{S}\cap \hat{R}$, then both $L\cup (\hat{S}\setminus\hat{R})$ and $L\cup (\hat{R}\setminus\hat{S})$ are independent sets. Moreover, as $L\subseteq V(G_{\mathcal{B}})\setminus V_0$, $L$ is  an IT of $\mathcal{B}_0\setminus\{V_0\}$. Hence, both $L\cup (\hat{S}\setminus\hat{R})$ and $L\cup (\hat{R}\setminus\hat{S})$ are ITs of $\mathcal{B}'_0$. Finally, since $D'\setminus V(C')= D\setminus V(C)$, we have $|D'\setminus V(C')|< \epsilon (|\mathcal{B}|-1)\leq\epsilon (|\mathcal{B'}|-1)$. 

\end{itemize}
This completes the proof of Lemma~\ref{domin_set_G'}.
\end{proof}

\section{Proofs: $S$--$R$ Path}\label{sec:proof:SRPath}
In this section, we show that \textsc{FindSR-Path} terminates and returns a polynomial-length $S$--$R$ path $P$ in $\textbf{RG}(G,\mathcal{V})$, provided that the condition in line $26$ of \textsc{GrowPath} is not satisfied at any point during the execution of \textsc{FindSR-Path}. We begin by bounding the number of layers in the alternating tree.

\begin{lemma}\phantomsection\label{fewlayers} Suppose $(\mu,U,\rho)$ is feasible for $(r,\epsilon)$. Then, throughout the execution of {\rm \textsc{GrowPath}}, the alternating tree maintained by the algorithm has at most $c\log m$ layers, where \[ c=\frac{1}{\log(1+\rho(1-\mu))}. \] 
\end{lemma} 

\begin{proof} Let $T=(L_0,\ldots,L_\ell)$ be the alternating tree at the beginning of an iteration of the while loop in line~22. Fix $i\in [\ell]$. By Lemma~\ref{no_collapse}, $ |Y_i|\ge (1-\mu)|X_i|,$ and, since {\rm \textsc{GrowPath}} did not terminate when $L_i$ was created, Lemma~\ref{domin_set_G} yields $ |X_i|> \rho |Y_{\le i-1}|.$ Combining these inequalities gives $ |Y_i| > \rho(1-\mu)|Y_{\le i-1}|. $ Consequently, \[ |Y_{\le i}| = |Y_{\le i-1}|+|Y_i| > \bigl(1+\rho(1-\mu)\bigr)|Y_{\le i-1}|. \] Note that, by $(\hyperref[a_2]{\text{a$_2$}})$, we may assume that $A\neq \emptyset$ which implies $|X_1|\geq 1$. Consequently, we have $|Y_{1}|\ge (1-\mu)|X_1|> 0$ and hence $|Y_{1}|\ge 1$. Repeated application of the above inequality then yields \[ |Y_{\le\ell}| > \bigl(1+\rho(1-\mu)\bigr)^\ell. \] On the other hand, every vertex of $Y_{\le\ell}$ belongs to $M$, so $ |Y_{\le\ell}|\le m.$ Therefore, $ m \ge \bigl(1+\rho(1-\mu)\bigr)^\ell,$ which is equivalent to \[ \ell \le \frac{\log m}{\log(1+\rho(1-\mu))} = c\log m. \] This proves Lemma~\ref{fewlayers}.
\end{proof}

We now turn to the signature vector and show that it makes strict progress in every non-terminating iteration.
\begin{lemma}\phantomsection\label{lexreduce}
The lexicographic value of the signature vector decreases across each
iteration of the loop in line $22$ of {\rm \textsc{GrowPath}}, unless the
algorithm terminates during that iteration.
\end{lemma}

\begin{proof}
Let $T=(L_0,\ldots,L_\ell)$ be the alternating tree at the beginning of
an iteration of the while loop in line $22$ of {\rm \textsc{GrowPath}},
and let
$s=(s_1,\ldots,s_{2\ell},\infty)$ be its signature vector.

We distinguish two cases.

\begin{case}
No collapse operation occurs, i.e.,
$|I_M(X_{\ell+1})|\le\mu|X_{\ell+1}|$.
\end{case}

The only modification to $T$ is that a new layer $L_{\ell+1}$ is
appended. Hence the new signature vector is obtained from $s$ by
appending two coordinates before the terminal $\infty$. Thus $s$ is a
proper prefix of the new signature vector, and therefore its
lexicographic value decreases.

\begin{case}
At least one collapse operation occurs, i.e.,
$|I_M(X_{\ell+1})|>\mu|X_{\ell+1}|$.
\end{case}

Let $t$ be the index of the layer satisfying
$|I_M(X_t)|>\mu|X_t|$ in the final iteration of the collapse loop (lines 30--51) of the current iteration of $\textsc{GrowPath}$. If $t=1$, then a vertex of $I_M(X_1)$ is added to $M$ and the
algorithm terminates. Thus, assume $t>1$. Let $\hat{T}$ denote the alternating tree returned by $\textsc{SuperposedBuild}$ in this final iteration when called with input tree $T'$.

\subcase{2.1.}
{\rm During this final collapse, \textsc{SuperposedBuild}} returns $T'$ as $\hat{T}$.

 Since $T'$ has $t-1$ layers, all signature coordinates preceding $s_{2t-3}$ remain unchanged. Moreover, since $X_{i}$'s are not modified in {\rm \textsc{GrowPath}} by lines $39-46$, we have $X'_{t-1}=X_{t-1}$, and hence $s_{2t-3}$ also remains unchanged. Therefore, all signature coordinates preceding $s_{2t-2}$ remains unchanged. Hence, the signature vector of $\hat{T}=T'$ is $'=(s'_1,\ldots,s'_{2t-2},\infty)$ such that $s'_i = s_i$ for all $1 \le i \le 2t-3$, and \[ s'_{2t-2} =  \left\lfloor \log_b \frac{\rho^{-(t-1)}}{(1-\mu)^{t-1}}|Y'_{t-1}| \right\rfloor. \]

Every vertex class containing a vertex of $I_M(X_t)$ also contains a vertex of $Y_{t-1}$. Since each vertex class contains at most $U$ vertices of $X_t$, at least $\frac{\mu}{U}|X_t|$ blocking vertices are removed from $Y_{t-1}$ by the replacements in lines $39-46$. Note that we perform these replacements one by one, and this procedure does not change the set of vertices that are candidates for replacement. More precisely, let $\mathcal{W}=\{w_1,\ldots,w_k\}$ be the set of all vertices $w\in Y_{t-1}$ such that $I_M(X_t)\cap\operatorname{cls}(w)\neq\emptyset$, and, for each $i\in[k]$, let $u_i\in I_M(X_t)\cap\operatorname{cls}(w_i)$. Suppose that $M_{\text{new}}=(M\setminus\{w_1\})\cup\{u_1\}$. Now, consider any $i\in[k]\setminus\{1\}$. We still have $w_i\in Y_{t-1}$. Moreover, based on \textsc{BuildLayer}, which produces $X_t$, every new vertex added to $X_t$ must be addable with respect to the previous vertices, and hence $X_t$ is an independent set. Therefore, $u_i\notin N[u_1]$, and since $u_i\in I_M(X_t)\subseteq I_{M\setminus\{w_1\}}(X_t)$, we conclude that $u_i\in I_{M_{\text{new}}}(X_t)$. Hence, for every $w_i\in\mathcal{W}\setminus\{w_1\}$, we still have $I_{M_{\text{new}}}(X_t)\cap\operatorname{cls}(w_i)\neq\emptyset$. Thus, performing the replacements one by one does not affect the remaining candidates for replacement, ensuring that at least $\frac{\mu}{U}|X_t|$ blocking vertices are removed from $Y_{t-1}$, as claimed. Therefore, as we assume that {\rm \textsc{GrowPath}} has not terminated, we have:

\begin{align*}
|Y'_{t-1}|
&\le
|Y_{t-1}|-\frac{\mu}{U}|X_t|&\\
&\leq |Y_{t-1}|-\frac{\mu\rho}{U}
|Y_{\leq t-1}|&\hspace{3em}({\text{Lemma~\ref{domin_set_G}}})&\\
&\le |Y_{t-1}|-\frac{\mu\rho}{U}
|Y_{t-1}|=
\left(1-\frac{\mu\rho}{U}\right)|Y_{t-1}|.&
\end{align*}

Hence,

\begin{align*}
\log_b \frac{\rho^{-(t-1)}}{(1-\mu)^{t-1}}|Y'_{t-1}|&< \log_b \frac{\rho^{-(t-1)}}{(1-\mu)^{t-1}}\left(1-\frac{\mu\rho}{U}\right)|Y_{t-1}| \\
&\le  \log_b \frac{\rho^{-(t-1)}}{(1-\mu)^{t-1}}|Y_{t-1}| + \log_b \left(1-\frac{\mu\rho}{U}\right)\\
&=  \log_b \frac{\rho^{-(t-1)}}{(1-\mu)^{t-1}}|Y_{t-1}| + \log_b b^{-1}=\log_b \frac{\rho^{-(t-1)}}{(1-\mu)^{t-1}}|Y_{t-1}|-1.
\end{align*}

which implies
$s'_{2t-2}<s_{2t-2}$.Thus, the lexicographic value decreases.

\subcase{2.2.}
{\rm \textsc{SuperposedBuild}} returns $\hat{T}$ distinct from $T'$ .

Let $q$ be the index of the last layer of $\hat{T}$. Then
\[
|X'_q|\ge(1+\mu)|X_q|.
\]
Since $(\mu,U,\rho)$ is feasible, by $(\hyperref[e_3]{\text{e$_3$}})$, we have $\frac{\mu\rho}{U-\mu\rho}<\mu$. So, $b=\frac{U}{U-\mu\rho}=1 +\frac{\mu\rho}{U-\mu\rho}< 1+\mu,$
and therefore,
\begin{align*}
\log_b \frac{\rho^{-q}}{(1-\mu)^{q-1}}|X'_q|&\ge \log_b \frac{\rho^{-q}}{(1-\mu)^{q-1}}(1+\mu)|X_{q}| \\
&=\log_b \frac{\rho^{-q}}{(1-\mu)^{q-1}}|X_{q}|+\log_b\left(1+\mu\right)\\
&>\log_b \frac{\rho^{-q}}{(1-\mu)^{q-1}}|X_{q}| +1.
\end{align*}
Hence $s'_{2q-1}< s_{2q-1}$. No layer preceding $L_q$ is modified, so all earlier signature
coordinates remain unchanged. Hence, the signature vector of $\hat{T}$ is $s'=(s'_1,\ldots,s'_{2q-1},s'_{2q},\infty)$ such that $s'_i = s_i$ for all $1 \le i \le 2q-2$, and \[ s'_{2q-1} =  -\left\lfloor \log_b \frac{\rho^{-q}}{(1-\mu)^{q-1}}|X'_q| \right\rfloor. \]

Therefore, the first changed coordinate is smaller, and the signature vector decreases lexicographically.
\end{proof}

To bound the number of possible signature vectors, we first establish a monotonicity property of their coordinates.

\begin{lemma}\phantomsection\label{abso_valu} The coordinates of the signature vector are non-decreasing in absolute value at the beginning of each iteration of the while loop in line $22$ of {\rm \textsc{GrowPath}}. 
\end{lemma} 

\begin{proof} Let $T=(L_0,\ldots,L_\ell)$ be the alternating tree with respect to the PIT $M$ at the beginning of an iteration of the while loop. It suffices to verify that every pair of consecutive coordinates of the signature vector is non-decreasing in absolute value. First consider the two coordinates corresponding to the same layer $L_i$, where $i\in[\ell]$. By Lemma~\ref{no_collapse}, \[ |Y_i|\ge (1-\mu)|X_i|, \] and therefore \[ \frac{\rho^{-i}}{(1-\mu)^{i-1}}|X_i| \le \frac{\rho^{-i}}{(1-\mu)^i}|Y_i|. \] Since $\lfloor\log_b(\cdot)\rfloor$ is monotone increasing, \[ |s_{2i-1}| = \left\lfloor \log_b\!\left(\frac{\rho^{-i}}{(1-\mu)^{i-1}}|X_i|\right) \right\rfloor \le \left\lfloor \log_b\!\left(\frac{\rho^{-i}}{(1-\mu)^i}|Y_i|\right) \right\rfloor = |s_{2i}|. \] Next consider two consecutive coordinates coming from adjacent layers, namely $s_{2i}$ and $s_{2i+1}$ for $i\in[\ell-1]_0$. By Lemma~\ref{domin_set_G}, \[ |X_{i+1}|\ge \rho |Y_i|, \] which implies \[ \frac{\rho^{-i}}{(1-\mu)^i}|Y_i| \le \frac{\rho^{-(i+1)}}{(1-\mu)^i}|X_{i+1}|. \] Applying the monotonicity of $\lfloor\log_b(\cdot)\rfloor$ again gives \[ |s_{2i}| = \left\lfloor \log_b\!\left(\frac{\rho^{-i}}{(1-\mu)^i}|Y_i|\right) \right\rfloor \le \left\lfloor \log_b\!\left(\frac{\rho^{-(i+1)}}{(1-\mu)^i}|X_{i+1}|\right) \right\rfloor = |s_{2i+1}|. \] Thus every consecutive pair of coordinates satisfies the desired inequality, and hence the coordinates of the signature vector are non-decreasing in absolute value. \end{proof}

Combining this monotonicity with the bound on the number of layers, we can now bound the total number of possible signature vectors.

\begin{lemma}\phantomsection\label{fewsigs} Let $T=(L_0,\ldots,L_\ell)$ and $M$ be the alternating tree and the partial independent transversal at the beginning of an iteration of the while loop in line~22 of {\rm \textsc{GrowPath}}. The collection of all possible signature vectors of $T$ has cardinality polynomial in $m$, with exponent depending only on $r$ and $\epsilon$. \end{lemma}
\begin{proof} Since we have at most $c\log m$ layers (Lemma~\ref{fewlayers},) and each layer of the alternating tree contributes exactly two entries to its signature vector, the total number of finite coordinates in the signature vector is at most $2c\log m$. Next, Lemma~\ref{abso_valu} implies that the coordinates of the signature vector are monotone with respect to their absolute values. Therefore, the last finite coordinate has the largest absolute value among all coordinates. By Definition~\ref{sig}, this coordinate equals $$ \left\lfloor \log_b\left( \frac{\rho^{-\ell}}{(1-\mu)^\ell}|Y_\ell| \right) \right\rfloor. $$ Since $\ell\le c\log m$ ( Lemma~\ref{fewlayers}) and $|Y_\ell|\le m$, every coordinate has absolute value at most \begin{align*} \log_b\left( \frac{\rho^{-c\log m}}{(1-\mu)^{c\log m}}m \right) &=\log_b m+\log_b\!\left(\rho^{-c\log m}\right) -\log_b\!\left((1-\mu)^{c\log m}\right)\\ &=\frac{\log m}{\log b} -c\log_b(\rho)\log m -c\log_b(1-\mu)\log m\\ &= \left( \frac1{\log b} -c\log_b(\rho) -c\log_b(1-\mu) \right)\log m. \end{align*} Define $$ R= \frac1{\log b} -c\log_b(\rho) -c\log_b(1-\mu). $$ Observe that $R$ depends only on $r$ and $\epsilon$, since the parameters $b$, $c$, $\mu$, $U$, and $\rho$ are determined solely by $r$ and $\epsilon$. Now consider an arbitrary signature vector $ s=(s_1,s_2,\ldots,s_{2\ell},\infty),$ and define the shifted vector $$ s^+= (s_1-1,s_2+2,\ldots,s_{2\ell-1}-(2\ell-1),s_{2\ell}+2\ell,\infty). $$ The largest finite coordinate of $s^+$ is bounded by $ R\log m+2\ell \le (R+2c)\log m. $ Furthermore, because the coordinates of $s$ are non-decreasing in absolute value with the prescribed alternating sign pattern, the coordinates of $s^+$ become strictly increasing in absolute value. Consequently, every vector $s^+$ uniquely determines a subset of $$ \{1,\ldots,\lfloor(R+2c)\log m\rfloor\}, $$ and distinct signature vectors produce distinct subsets. Therefore, the number of possible vectors $s^+$, and hence the number of possible signature vectors, is at most $ 2^{(R+2c)\log m},$ which is a polynomial in $m$ whose degree depends only on $r$ and $\epsilon$. 
\end{proof}
Now, we are ready to prove the correctness of Theorem~$\hyperref[thm:mainprime]{1.4}$, leaving aside the polynomial running-time analysis, which will be established in Section~\ref{sec:proof:runningtime}.

\textbf{Correctness of Theorem~$\hyperref[thm:mainprime]{1.4}$:}
Let $r\in \poi$, $\epsilon >0$, and $G'$ be a $r$-claw-free graph with a vertex partition $\mathcal{V'}=(V'_1, V'_2,\dots,V'_m)$, and two independent transversals $S,R\in \mathcal{I}(G',\mathcal{V'})$ are given. We provide $(G',\mathcal{V}',S,R)$ as input to \textsc{FindSR-Path}. During any call to \textsc{GrowPath}, if either the condition in line 19 or that in line 26 is satisfied, then the second outcome of Theorem~\ref{thm:main} immediately follows by $(\hyperref[a_2]{\text{a$_2$}})$ or Lemma~\ref{domin_set_G'}, respectively. Thus, we may assume that $A\neq\emptyset$ and 
$|X_{\ell+1}|>\rho|Y_{\leq\ell}|$ throughout the execution of \textsc{GrowPath}. Consider one iteration of the loop in line 4 of \textsc{FindSR-Path}. Under this assumption, $\textsc{GrowPath}$ eventually terminates with $x=0$ at either line 17 or line 37, since each iteration of the while loop in line 22 strictly decreases the lexicographic value of the signature vector and there are only finitely many possible signature vectors (Lemmas~\ref{lexreduce} and ~\ref{fewsigs}). Consequently, $|\tilde{S}\cap \tilde{R}|$ must be augmented either in lines 13--14 or in lines 33--34 of \textsc{GrowPath}. In particular, the size of $\tilde{S}\cap \tilde{R}$, which is at least $|\hat{S}\cap \hat{R}|$, increases by one (see $(\hyperref[a_1]{\text{a$_1$}})$ for details of lines 13--14). Moreover, the changes from $\hat{S}$ to $\tilde{S}$ and from $\hat{R}$ to $\tilde{R}$ each constitute a valid single reconfiguration step. Hence, $\hat{S}\tilde{S}$ and $\hat{R}\tilde{R}$ are edges of $\textbf{RG}(G',\mathcal{V}')$. Consider now one iteration of the for-loop in line 39 of \textsc{GrowPath}. Let $M_\text{new}$, $\title{S}_\text{new}$, and $\title{R}_\text{new}$ be the PIT $M$ and two ITs $\title{S},\title{R}$ at the end of the current iteration, respectively. Similarly, let $M_\text{old}$ $\title{S}_\text{old}$, and $\title{R}_\text{old}$ denote their values at the end of the preceding iteration. Since $M_\text{old}$ and $M_\text{new}$ differ by exactly one element, the corresponding changes from $\title{S}_\text{old}$ to $\title{S}_\text{new}$ and from $\title{R}_\text{old}$ to $\title{R}_\text{new}$ each constitute a valid single reconfiguration step. Hence, $\title{S}_\text{new}$$\title{S}_\text{old}$ and $\title{R}_\text{new}$$\title{R}_\text{old}$ are two edges of $\textbf{RG}(G', \mathcal{V}')$. Therefore, in either case (line 17 or line 37), the output $P_{\hat{S}}$ and $P_{\hat{R}}$ are two paths starting from $\hat{S}$ and $\hat{R}$, respectively, such that $|\text{End($P_{\hat{S}}$)}\ \Delta \ \text{End($P_{\hat{R}}$)}| < |\hat{S}\Delta\hat{R}|$. Regarding the length of this $\hat{S}$--$\hat{R}$ path, by Lemma~\ref{fewsigs} and Lemmas~\ref{lexreduce}, we know that the while loop in line~22 of \textsc{GrowPath} repeats at most $P(m)$ times, where $P(m)$ is a polynomial in $m$ whose degree depends only on $r$ and $\epsilon$. During each iteration of the while loop in line~22, the while loop in line~30 repeats at most $\ell=c\log m$ times (see Lemma~\ref{fewlayers} for the definition of $c$), since $\ell$ decreases in each iteration of the latter loop. Moreover, during each iteration of the while loop in line~30, the for loop in line~39 repeats at most $m$ times, since $Y_{\ell-1}\subseteq M$. Therefore, after a call to \textsc{GrowPath}, the lengths of the paths $P_{\hat{S}}$ and $P_{\hat{R}}$ are $O(cm\log mP(m))$, which is polynomial in $m$. Since these two paths meet after at most $m$ iterations of the loop in line~4 of \textsc{FindSR-Path}, the total length of the resulting path is $O(cm^2\log mP(m))$, which is polynomial in $m$ when $r$ and $\epsilon$ are regarded as fixed constants.

\section{Proofs: Running Time}\label{sec:proof:runningtime}
As the final part of the proof, we show that the termination of \textsc{FindSR-Path}, established in Section~\ref{sec:proof:SRPath}, occurs after a polynomial number of iterations, where the polynomial is in $m$ and its degree depends only on $r$ and $\epsilon$, which are regarded as fixed constants.

 Recall that $(\mu,U,\rho)=(\frac{\epsilon}{10r},\frac{10r}{\epsilon},\frac{\epsilon}{10r})$ is feasible for $(r,\epsilon)$ when $r\geq 2$ and $0<\epsilon<1$. Under this choice of parameters, the following Lemma regarding the parameter $\mu$ helps us to obtain a closed form for the running time of the algorithm.

\begin{lemma}\label{e456}
Let $\mu=\frac{\epsilon}{10r}$ where $r\geq 2$ and $0<\epsilon<1$. Then $\mu$ satisfies the following three inequalities;
\begin{enumerate} 
\item[\textnormal{e\textsubscript{4}:}]\label{e_4} $\log\left(\frac{1}{1-\mu^3}\right)\geq \mu^3,$
\vspace{0.5em}
\item[\textnormal{e\textsubscript{5}:}]\label{e_5} $\log(1+\mu-\mu^2)\geq \frac{\mu}{4},$
\vspace{0.5em}
\item[\textnormal{e\textsubscript{6}:}]\label{e_6} $\log\left(\frac{1+\mu-\mu^2}{\mu-\mu^2}\right) +2\log\left(\frac{1}{1-\mu^3}\right) \leq 2\log\left(\frac{3}{\mu}\right).$ 

\end{enumerate}
\end{lemma}
\begin{proof} 
Since \(r\geq 2\) and \(0<\epsilon<1\), we have $ 0<\mu<\frac12.$ We first prove~\hyperref[e_4]{\text{e$_4$}}. Recall that $-\log(1-t)\geq t$ for every $0\leq t<1$. Applying this inequality with \(t=\mu^3\)  proves~\hyperref[e_4]{\text{e$_4$}}.

To prove~\hyperref[e_5]{\text{e$_5$}}, set  $t=\mu-\mu^2=\mu(1-\mu)$. Since $0<\mu<\frac12$, we have  $\frac{\mu}{2}<t<1$. Moreover, we have $ \log(1+t)\geq \frac{t}{2}$ for every $0\leq t\leq 1$. It follows that \[ \log(1+\mu-\mu^2) =\log(1+t) \geq \frac{t}{2} \geq \frac{\mu}{4}, \] and hence~\hyperref[e_5]{\text{e$_5$}} holds. Finally, we prove~\hyperref[e_6]{\text{e$_6$}}. Since $ 1-\mu\geq \frac12$ and $1+\mu-\mu^2\leq \frac32$, we obtain \[ \frac{1+\mu-\mu^2}{\mu-\mu^2} = \frac{1+\mu-\mu^2}{\mu(1-\mu)} \leq \frac{3}{\mu}. \] Therefore, \[ \log\left(\frac{1+\mu-\mu^2}{\mu-\mu^2}\right) \leq \log\left(\frac{3}{\mu}\right). \] Furthermore, since $\mu<\frac12$, we have $\mu^3<\frac18$, and thus  $\frac{1}{1-\mu^3}\leq \frac87$. Also, $\left(\frac87\right)^2 =\frac{64}{49} <2 <\frac{3}{\mu}$. Consequently, \[ 2\log\left(\frac{1}{1-\mu^3}\right) \leq 2\log\left(\frac87\right) \leq \log\left(\frac{3}{\mu}\right). \] Adding the previous two inequalities gives \[ \log\left(\frac{1+\mu-\mu^2}{\mu-\mu^2}\right) +2\log\left(\frac{1}{1-\mu^3}\right) \leq 2\log\left(\frac{3}{\mu}\right), \] which proves~\hyperref[e_6]{\text{e$_6$}}. 
\end{proof}

\textbf{Running Time:}
Now, we are ready to analyze the running time of the algorithm by considering each subroutine separately. It is clear that {\rm \textsc{BuildLayer}} terminates in time $O(n^2)$, since testing whether a given vertex is addable can be checked in $O(n)$ time. Since the loop in line~3 is executed at most $\ell$ times, {\rm \textsc{SuperposedBuild}} runs in time $O(\ell n^2)$.
For {\rm \textsc{GrowPath}}, one can easily check that lines~2--11 can be executed in $O(n^2)$ time. By \cref{lexreduce}, the lexicographic value of the signature vector decreases at each iteration of the loop in line~22. Moreover, by \cref{fewsigs}, the number of possible signature vectors is at most \[ 2^{(R+2c)\log m}=m^{R+2c}, \] where 
\begin{align*} R=\frac{1}{\log b}-c\log_b(\rho)-c\log_b(1-\mu),\text{ and } \ c=\frac{1}{\log(1+\mu-\mu^2)}. 
\end{align*}
This implies that the loop in line 22 is executed at most $m^{R+2c}$ times. Since all other steps within the loop in line~22 can be performed in time $O(\ell^2n^2)$, the loop in line~22, and consequently {\rm \textsc{GrowPath}}, terminates in time \[ O\left(m^{R+2c}\ell^2n^2\right). \] Finally, for {\rm \textsc{FindSR-Path}}, all steps of the loop in line~4, except for the call in line~5, can be performed in time $O(\ell n^2)$. Since the loop in line~4 is executed at most $m$ times and calls {\rm \textsc{GrowPath}} in line~5 each time, the total running time of {\rm \textsc{FindSR-Path}} is \[ O\left(m^{R+2c+1}\ell^2n^2\right). \]

Now, we have

\begin{align*}
R&= \frac1{\log b} -c\log_b(\rho) -c\log_b(1-\mu) \\
&=\frac1{\log b} -c\frac{\log\rho}{\log b} -c\frac{\log(1-\mu)}{\log b} \\
&=\frac{\log(1+\rho(1-\mu))-\log\rho-\log(1-\mu)}{\log (\frac{1}{1-\mu^3}).\log(1+\rho(1-\mu))}\\
&=\frac{\log(1+\mu-\mu^2)-\log\mu-\log(1-\mu)}{\log (\frac{1}{1-\mu^3}).\log(1+\mu-\mu^2)}\\
&=\frac{\log(\frac{1+\mu-\mu^2}{\mu(1-\mu)})}{\log (\frac{1}{1-\mu^3}).\log(1+\mu-\mu^2)},
\end{align*}

and so, 

\begin{align*}
R+2c&=\frac{\log(\frac{1+\mu-\mu^2}{\mu(1-\mu)})}{\log (\frac{1}{1-\mu^3}).\log(1+\mu-\mu^2)}+\frac{2}{\log(1+\mu-\mu^2)}\\
&=\frac{\log(\frac{1+\mu-\mu^2}{\mu(1-\mu)})+2\log (\frac{1}{1-\mu^3})}{\log (\frac{1}{1-\mu^3}).\log(1+\mu-\mu^2)}.
\end{align*}
\\
Thus, by Lemma~\ref{e456}, we obtain

\[
R + 2c\leq \frac{8}{\mu^4}\log(\frac{3}{\mu}),
\]
and since $\frac{1}{\mu}>2$, we conclude that,

\begin{align*}    
R + 2c+1+2\log c&\leq \frac{8}{\mu^4}\log(\frac{3}{\mu}) +1 +2\log (\frac{4}{\mu})\\
&\leq \frac{8}{\mu^4}\log(\frac{3}{\mu}) +3\log (\frac{4}{\mu})\\
&\leq \frac{8}{\mu^4}\log(\frac{3}{\mu}) +3(\log_{3}4)\log (\frac{3}{\mu})\\
&\leq \frac{8}{\mu^4}\log(\frac{3}{\mu}) +\frac{1}{\mu^2}(\frac{1}{\mu^2})\log (\frac{3}{\mu})\\
&\leq \frac{9}{\mu^4}\log(\frac{3}{\mu}).
\end{align*}

Recall that, by Lemma~\ref{fewlayers}, we have $\ell=c\log m$. Hence, the total running time is
\begin{align*}
O\left(m^{R+2c+1}\ell^2n^2\right)&=O\left(m^{R+2c+1}(c^2\log^2 m)n^2\right)\\
&=O\left(m^{R+2c+1+2\log c}(\log^2 m)n^2\right)\\
&=O\left(m^{\frac{9}{\mu^4}\log(\frac{3}{\mu})}(\log^2 m)n^2\right).\\
\end{align*}
\\

\section{Application: Polynomial Path Between Two Hypergraph Matchings }\label{sec:Application}
In this section, we apply Theorem~\ref{thm:main} to obtain a polynomial-time algorithm for reconfiguring perfect matchings in uniform bipartite hypergraphs. Throughout this section, we use the term \textit{hyperedges} for edges of a hypergraph, to avoid confusion with the edges in an auxiliary graph that will be defined.
An \emph{$r$-uniform bipartite hypergraph} is a hypergraph $H=(A,B,E)$ whose vertex set is divided into two parts $A$ and $B$, where every hyperedge contains exactly one vertex from $A$ and $r-1$ vertices from $B$. A \emph{perfect matching} of $H$ is a family $\mathcal{M}\subseteq E$ of pairwise disjoint hyperedges that covers all vertices of $A$. Since each hyperedge contains precisely one vertex of $A$, this is equivalent to requiring  $|\mathcal{M}|=|A|$. Given $S\subseteq A$, we define the \textit{link} $E_S=\{e\in E:e\cap S\neq\emptyset\}$ as the collection of hyperedges that meet $S$. Thus every hyperedge in $E_S$ intersects $S$ in exactly one vertex. Finally, for any family of hyperedges $F\subseteq E$, we define its \textit{cover number} $\tau_B(F)$ in $B$ by \[ \tau_B(F) = \min\bigl\{|T|:T\subseteq B \text{ and } e\cap T\neq\emptyset \text{ for every }e\in F\bigr\}. \]

Haxell~\cite{Haxell_Match_HyperG} showed that an $r$-uniform bipartite hypergraph $H$ has a perfect matching if $\tau_B(E_S)>(2r-3)(|S|-1)$ for every $S\subseteq A$, and proved that this bound is best possible. Our result strengthens this existence theorem by showing that, under a slightly stronger condition, essentially obtained by adding an $\epsilon$ to the bound, we can also efficiently reconfigure between any two perfect matchings.

Given $H=(A,B,E)$, we associate with $H$ an auxiliary graph $G^H$ whose vertex set is  $V(G^H)=E$. Two distinct vertices $e,f\in V(G^H)$ are joined by an edge in $G^H$ precisely when the corresponding hyperedges have a common vertex in $B$, namely, \[ ef\in E(G^H) \quad\Longleftrightarrow\quad e\cap f\cap B\neq\emptyset. \] We further consider $G^H$ with a natural partition of its vertex set. For every $a\in A$, define $V_a=\{e\in E:a\in e\}.$ Since every hyperedge of $H$ contains exactly one vertex of $A$, the collection  $\mathcal{A}=\{V_a:a\in A\}$  forms a partition of $V(G^H)$. Thus, two hyperedges belong to the same class of $\mathcal{V}$ exactly when they contain the same vertex of $A$. Under this correspondence, perfect matchings in $H$ are precisely independent transversals of $(G^H,\mathcal{A})$. Moreover, observe that $G^H$ is $r$-claw-free (with respect to any vertex partition, in fact).

\begin{lemma}\label{cover}
Let $\mathcal{B}$ be a set of partition classes of $G^H$, and let $C$ be a
biconstellation for $\mathcal{B}$. Then 
$|\bigcup_{e\in V(C)}e\cap B|\leq(2r-3)|{\mathcal{B}}|$.
\end{lemma}
\begin{proof}[Proof:]
 Let $F$ denote the balanced bipartite graph of $C$ with each part of size $t$. Consider any star component $K$ of $C$. Then $K$ corresponds to a
set of hyperedges of $H$, consisting of the centre $e_K$ of the star $K$
and a nonempty set 
$L_K$ of leaves, all of which intersect $e_K$ in $B$. Hence the total
number of vertices of $B$ contained in $\{e_K\}\cup L_K$ is at most 
$(r-1)+(r-2)|L_K|$. By definition of $C$, $\bigcup_KL_K$ is an IT of $|\mathcal{B}|-t$
classes of $\mathcal{B}$. So, $\sum\limits_{K\in C\setminus F}|L_K|=|\mathcal{B}|-t$ and $C\setminus F$ has 
at most $|\mathcal{B}|-t$ star components. Therefore 
\begin{align*}
\left|\bigcup_{e\in V(C)}e\cap B\right|&\leq\left|\bigcup_{e\in V(F)}e\cap B\right|+\left|\bigcup_{e\in V(C\setminus F)}e\cap B\right|\\
&\leq t((r-1)+(r-2)) +\sum_{K}\left(r-1+(r-2)|L_K|\right)\\
&=(r-1)t+(r-2)t+ (r-1)(|\mathcal{B}|-t) + (r-2)(|\mathcal{B}|-t) \\
&=(r-1)|\mathcal{B}| + (r-2)|\mathcal{B}|=(2r-3)|\mathcal{B}|,
\end{align*}
where the sum is over all components $K$ of $C\setminus F$. The bound on the size of $\bigcup_{e\in V(F)} e\cap B$ follows from the fact that $F$ has no isolated vertices, so every vertex of $F$ has a neighbour in $F$ (see $(\hyperref[a_1]{\text{a$_1$}})$ in Lemma~\ref{domin_set_G'}).

\end{proof}

\begin{theorem}\label{thm:hyper}
    Let $r\in \poi$ and $\epsilon >0$ be given. Then there is an algorithm $A_{r,\epsilon}$ that takes as an input an $r$-uniform bipartite hypergraph $H=(A,B,E)$ and two perfect matchings $\mathcal{M}_1$ and $\mathcal{M}_2$  of $H$, and finds in polynomial time either 
   
    \begin{itemize}
        \item[-] a polynomial-length $\mathcal{M}_1$--$\mathcal{M}_2$ path $P$ in $\textbf{RG}(G^H, \mathcal{A})$, or,
        \item[-] a subset $S\subseteq A$ where $\tau_B(E_S)\leq(2r-3+\epsilon)|S|$.
    \end{itemize} 
\end{theorem}
\begin{proof}
Define $\epsilon'=\epsilon/(r-1)$ and apply Theorem~\ref{thm:main} with $r$ and $\epsilon'$ to obtain a polynomial-time algorithm $\mathcal{A}_{r,\epsilon'}$. Since the graph $G^H$ is $r$-claw-free for every $r$-uniform bipartite hypergraph $H=(A,B,E)$ with respect to any partition, we can feed $G^H$ to $\mathcal{A}_{r,\epsilon'}$ and obtain either

\begin{enumerate} 
\item[$\text{o}_1$:]\label{} a polynomial-length $\mathcal{M}_1$--$\mathcal{M}_2$ path $P$ in $\textbf{RG}(G^{H}, \mathcal{A})$, or 

\item[$\text{o}_2$:]\label{o_2} a set of partition classes $\mathcal{B}\subseteq \mathcal{A}$ and a vertex set $D\subseteq V(G^H)$ such that $|D|\leq (2+\epsilon')|\mathcal{B}|$ and $D$ dominates $G^H_{\mathcal{B}}$. Furthermore, there exists a biconstellation $C$ for some $\mathcal{B}_0\supseteq \mathcal{B}$ where  $V(C)\subseteq D$ and $|D\setminus V(C)|< \epsilon'(|\mathcal{B}|-1)$.
\end{enumerate}

If $\text{o}_1$ is the outcome, then we have found the polynomial path between the perfect matchings $\mathcal{M}_1$ and $\mathcal{M}_2$ that we seek. Suppose $\text{o}_2$ holds and let $F$ be the balanced bipartite graph of $C$ with parts of size $t$. Let $S
\subseteq A$ be the the set of indices of partition classes of  $\mathcal{B}$. Then $D$ is a set of hyperedges of $H$ such that every hyperedge of $E_{S}$ intersects $\bigcup_{e\in D}e\cap B$. To estimate $\left|\bigcup_{e\in D}e\cap B\right|$, we first estimate $|\bigcup_{e\in V(C)}e\cap B|$. Note that for each star component $K$ of the biconstellation $C$, the number of vertices of $B$ contained in $\{e_K\}\cup L_K$ is at most $(r-1)+(r-2)|L_K|=1+(r-2)|V(K)|$, as in the proof of Lemma~\ref{cover}, and $|V(K)|\geq 2$. Hence, the number of star components of $C\setminus F$ is at most $(|V(C)|-2t)/2$.  Since $|D|<(2+\epsilon')(|{\mathcal B}|)$ we find 

\begin{align*} 
\left|\bigcup_{e\in V(C)}e\cap B\right|&\leq\left|\bigcup_{e\in V(F)}e\cap B\right|+\left|\bigcup_{e\in V(C\setminus F)}e\cap B\right|\\
&\leq t((r-1)+(r-2)) + \sum_K\left(1+(r-2)|V(K)|\right)\\ 
&\leq t((r-1)+(r-2))+ (|V(C)|
-2t)/2+(r-2)(|V(C)|-2t)\\
&= (2r-3)(2t)/2+ (2r-3)(|V(C)|
-2t)/2\\ &=(2r-3)|V(C)|/2\\ &=(2r-3)(|D|-|D\setminus V(C)|)/2\\
&<(2r-3)(1+\epsilon'/2)|{\mathcal B}|-(r-3/2)|D\setminus V(C)|, \end{align*} where the sum is over all components $K$ of $C\setminus F$. Since we have obviously  $|\bigcup_{e\in D\setminus V(C)}e\cap B|\leq (r-1)|D\setminus V(C)|$, we conclude that
\begin{align*}
\left|\bigcup_{e\in D}e\cap B\right |&\leq\left|\bigcup_{e\in V(D\setminus C)}e\cap B\right|+\left|\bigcup_{e\in V(C)}e\cap B\right|\\&<(r-1)|D\setminus V(C)|+(2r-3)(1+\epsilon'/2)|{\mathcal B}|-(r-3/2)|D\setminus V(C)|\\ &= |D\setminus V(C)|/2+(2r-3+\epsilon'(r-3/2))|\mathcal{B}|\\&<(\epsilon'/2)|\mathcal{B}|+(2r-3+\epsilon'(r-3/2))|\mathcal{B}| &(\hyperref[o_2]{\text{o$_2$}})\\
&\leq(2r-3+\epsilon'(r-1))|\mathcal{B}|= (2r-3+\epsilon)|\mathcal{B}|.
\end{align*}

Therefore, $\tau_B(E_{S})\leq(2r-3+\epsilon)|\mathcal{B}|=(2r-3+\epsilon)|S|$ which implies that $S$ is the desired set in the Theorem~\ref{thm:hyper}. 
\end{proof}

\section{Conclusion and Open Problems}\label{discussion} \subsection{Conclusion} In this paper, we developed a polynomial-time algorithm (Theorem~$\hyperref[thm:mainprime]{1.4}$) for the reconfiguration of independent transversals in $r$-claw-free graphs, for fixed $r$ and $\epsilon>0$. Given two independent transversals, our algorithm either constructs a polynomial-length reconfiguration path between them or identifies a set of partition classes whose induced subgraph admits a small strong dominating set. Moreover, the dominating set contains a biconstellation satisfying additional structural properties. Our result provides an algorithmic counterpart to recent connectivity theorems for independent transversals and establishes polynomial bounds on the diameter of the corresponding reconfiguration graph under suitable domination conditions. We also applied our approach to the reconfiguration of perfect matchings in bipartite uniform hypergraphs (Theorem~\ref{thm:hyper}).

\subsection{Open Problems} Our results raise several natural questions concerning the running times of our algorithms, their parameterized complexity, and the possibility of extending them to more general graphs. Throughout the following discussion, $n$ denotes the input size, and $f$ denotes an arbitrary function of its specified parameters. 

\paragraph{\textbf{Improving the dependence on $r$ and $\epsilon$:}} The running times of our algorithms are polynomial for fixed $r$ and $\epsilon$, but the exponents of these polynomials have a fourth-degree dependence on these parameters. A natural question is whether this dependence can be improved. In particular, we ask the following. 

\begin{problem} Can we obtain algorithms whose running-time exponents depend linearly on $r$ and $1/\epsilon$, rather than having a fourth-degree dependence?
\end{problem}

\paragraph{\textbf{Fixed-parameter tractability:}} The FPT aspects of our algorithm are similar to those of the algorithmic results of Graf and Haxell in ~\cite{GrafHaxell}{} (see also ~\cite{Graf_Harris_Haxell}): both algorithms are fixed-parameter tractable (FPT) when parameterized by $m$, $r$, and $1/\epsilon$. 
As an example of the applications of independent transversal results in parameterized complexity, Bessy et al.~\cite{bessy} used~\cite{GrafHaxell} to develop the rainbow matching technique, with which they obtained FPT algorithms and polynomial kernels for packing problems and their hitting counterparts.

Since $r$ and $1/\epsilon$ appear in the exponent of $m$ in the running time of our algorithm, it belongs to XP when parameterized by $r$ and $1/\epsilon$ alone. It is natural to ask whether this could be improved, by eliminating the dependence on $r$ and $1/\epsilon$ in the exponent of the running time, and thus obtaining an FPT reconfiguration algorithm with respect to these two parameters.

\begin{problem} Can we obtain an algorithm for the IT reconfiguration problem with running time $f(r,1/\epsilon)(n)^{O(1)}$, where the exponent of $n$ is independent of $r$ and $\epsilon$?
\end{problem} 

As a first step toward answering this question, one could consider parameterization by the maximum degree $\Delta$ instead of $r$. Since every graph of maximum degree $\Delta$ is $(\Delta+1)$-claw-free, this gives a more restricted setting in which to investigate fixed-parameter tractability. 


\paragraph{\textbf{Extending the algorithms to general graphs:}} The $r$-claw-free assumption is needed in our argument to bound the error term (see Claim~\ref{userclawfree} in Lemma~\ref{domin_set_G}): we may have to add to the dominating set given by $\bigcup_{i}(X_i\cup Y_i)$ a few more vertices that (together with their neighbors in $M$) were not added to sets $X_i,Y_i$ in line 6 of {\rm \textsc{SuperposedBuild}} because the condition in line 5 was not satisfied. The number of these vertices is small, but without the $r$-claw-free condition we could not guarantee that they have altogether few neighbors in $M$. However, this assumption does not seem essential, so it would be interesting to investigate whether it could be relaxed or removed entirely. In particular, Graf, Harris, and Haxell~\cite{Graf_Harris_Haxell} showed that, for every fixed $\epsilon>0$, an independent transversal can be found by a randomized polynomial-time algorithm when every vertex class has size at least $(2+\epsilon)\Delta$, i.e. the $r$-claw-free condition can be avoided in this case. This motivates the question of whether a similar result could be obtained for the reconfiguration version of the problem.

\begin{problem} For every fixed $\epsilon>0$, can we find a polynomial-length reconfiguration path between any two independent transversals in polynomial time, provided that every vertex class has size at least $(2+\epsilon)\Delta$? 
\end{problem}

More generally, we could ask whether the $r$-claw-free assumption can be removed from our algorithms, without imposing any lower bound on the sizes of the vertex classes, while preserving the guarantee of either finding a reconfiguration path or a small dominating set, with a running-time exponent depending only on $\epsilon$.

\subsection{Declaration on the Use of Artificial Intelligence} Artificial intelligence tools (ChatGPT (GPT-6)) were used solely for grammatical corrections, language editing, and improving the clarity and presentation of the manuscript. All mathematical results, including the theorems, proofs, algorithms, and mathematical arguments, were developed independently by the authors without the use of artificial intelligence.

\bibliographystyle{abbrv}
\bibliography{Ref}

\end{document}